\documentclass[twocolumn]{autart2}    % Enable this line and disable the 
\usepackage{graphicx}          % Include this line if your 
\usepackage{cite}
\usepackage{amsmath,amssymb,amsfonts}
\usepackage{algorithm,algorithmicx}
\usepackage{epstopdf}
\usepackage{algorithm}
\usepackage{algpseudocode}
\usepackage{mathrsfs}%
\usepackage{float}
\usepackage{xcolor}%
\usepackage{manyfoot}%
\usepackage{booktabs}%
\usepackage{listings}%
\usepackage{makecell}
\usepackage{booktabs,array,tabularx}
\usepackage{capt-of}
\newsavebox{\dimensionbox}
\newsavebox{\tablebox}

\begin{document}
\newtheorem{theorem}{Theorem}
\newtheorem{lemma}{Lemma}
\newtheorem{definition}{Definition}
\newtheorem{remark}{Remark}
\newtheorem{corollary}{Corollary}
\newtheorem{example}{Example}
\newtheorem{assumption}{Assumption}
\def\dom{{\rm dom}}
\newfloat{framework}{t}{lof}
\floatname{framework}{Framework}
\newcommand{\proof}{\emph{Proof: }}

\begin{frontmatter}

\title{Zeroth-order Algorithms for Distributed Online Optimization against Byzantine Edge Attacks\thanksref{footnoteinfo}} % Title, preferably not more 

\thanks[footnoteinfo]{Funding: This work was supported by the National Natural Science Foundation of China under Grants No. 72531001, No. 72201008 and No. 72131001. We also acknowledge the support of computing resources provided by the High-performance Computing Platform of Peking University.}

\author[pku]{Yuhang Liu}\ead{yuhang.liu@pku.edu.cn},    % Add the 
\author[pku]{Yujie Tang}\ead{yujietang@pku.edu.cn},               % e-mail address 
\author[pku,cor1]{Wenjun Mei}\ead{mei@pku.edu.cn}  % (ead) as shown
\corauth[cor1]{Corresponding author}
\address[pku]{Department of Control and Systems Engineering, Peking University, Beijing 100871, China}  % Please supply                                              
%\address[Rome]{Department of Control and Systems Engineering, Peking University, Beijing 100871, China}             % full addresses
%\address[Baiae]{Department of Control and Systems Engineering, Peking University, Beijing 100871, China}        % here.

                             % keyword list or with the 
                                          % help of the Automatica 
                                          % keyword wizard

\begin{abstract}                          % Abstract of not more than 200 words.

Byzantine attacks pose a fundamental challenge to distributed optimization, especially when they occur on communication edges, where truthful and corrupted messages coexist and attacked links are difficult to identify. While most existing studies focus on time-invariant objectives, many practical systems operate in dynamic environments and require distributed online optimization under such attacks. A central question is whether sublinear dynamic regret, a benchmark performance in online optimization, can still be achieved under Byzantine edge attacks. This paper answers this question affirmatively. 
%We directly study distributed online optimization %{in a challenging setting} 
%with time-varying directed networks, unknown and time-varying Byzantine edge attacks, and only zeroth-order feedback. 
{We directly study distributed online optimization with time-varying directed networks, where the identities of adversarial edges are unknown to the agents and may change over time, and only zeroth-order feedback is available. } 
Two resilient algorithms are proposed based on randomized and deterministic zeroth-order schemes. The former is more query-efficient and requires weaker smoothness conditions, while the latter provides deterministic guarantees and better performance in simulations.  {Both algorithms achieve consensus and, under suitable minimizer-path variation and algorithmic parameter conditions, achieve sublinear dynamic regret with respect to the original unweighted sum of local losses, without imposing stochastic assumptions for losses.} Numerical results validate the theoretical findings. Overall, the results show that reliable distributed online optimization remains possible even under complex Byzantine edge attacks and limited feedback.

%Byzantine attacks pose a fundamental challenge to distributed optimization, especially when they occur on communication edges, where truthful and corrupted messages coexist and attacked links are difficult to identify. While most existing studies focus on time-invariant objectives, many practical systems operate in dynamic environments and require distributed online optimization under such attacks. A central question is whether sublinear dynamic regret, a benchmark performance in online optimization, can still be achieved under Byzantine edge attacks. This paper answers this question affirmatively. We directly study distributed online optimization with time-varying directed networks, where the identities of adversarial edges are unknown to the agents and may change over time, and only zeroth-order feedback is available.  Two resilient algorithms are proposed based on randomized and deterministic zeroth-order schemes. The former is more query-efficient and requires weaker smoothness conditions, while the latter provides deterministic guarantees and better performance in simulations. Both algorithms achieve consensus and, under suitable minimizer-path variation and algorithmic parameter conditions, achieve sublinear dynamic regret with respect to the original unweighted sum of local losses, without imposing stochastic assumptions for losses. Numerical results validate the theoretical findings. Overall, the results show that reliable distributed online optimization remains possible even under complex Byzantine edge attacks and limited feedback.
\end{abstract}

\end{frontmatter}

\section{Introduction}
\label{sec:introduction}
\subsection{Background and motivation}

Distributed optimization has become a fundamental paradigm for large-scale multi-agent systems, particularly in cyber--physical systems, where agents cooperatively minimize a global objective through local information exchange~\cite{2009Nedic,Lobel2010,2014Chang}. However, cyber-attacks have emerged as one of the most critical threats in such systems, severely degrading the performance of distributed algorithms and even leading to system damage or paralysis~\cite{2019Zhao,2015Sundaram}. 
Among them, Byzantine attacks are particularly challenging due to their ability to generate arbitrary and inconsistent information~\cite{2019Sundaram}. 
{Byzantine attacks target either nodes or edges. In a node attack, an adversary controls a node and can send arbitrary messages to its
neighbors~\cite{2019Sundaram}. In an edge attack, an adversary can alter messages in transit on individual links~\cite{2023Xu}. Edge attacks pose distinct challenges for detection and mitigation because messages sent by honest nodes may be corrupted on specific links, while receivers typically do not know which links are compromised.
%Compared with node attacks, edge attacks are more subtle, as a corrupted received message does not imply that its sender is faulty, while the compromised edges are typically unknown. 
}

While most existing Byzantine-resilient optimization methods focus on time-invariant objectives, many modern cyber--physical systems, such as power systems and distributed target tracking, involve time-varying objectives~\cite{6702462,Aeron2008Efficient}. The study of distributed online convex optimization (DOCO)~\cite{2013Yan} integrates distributed optimization with online learning~\cite{2012Shalevshwartz,2016Hazan} to track time-varying optimal solutions using only past observations. A standard performance criterion in this setting is the achievement of sublinear regret. However,  whether such performance guarantees can still be ensured in the presence of Byzantine edge attacks remains unexplored.

%\begin{figure}[t]
%\begin{center}
%\includegraphics[scale=0.40]{figures/two-attack.pdf}    
%\caption{{Node attack and edge attack.}} \label{fig:node&edge}%%图片的label需要放在\caption之后（\caption 后才会更新图片编号） 
%\end{center}
%\end{figure}

{Meanwhile, many distributed optimization algorithms rely on gradient or subgradient information, which may be unavailable or computationally prohibitive to obtain in black-box settings and applications such as online advertising and graphical model inference~\cite{2016Hazan,2008Wainwright}. 
This motivates the use of zeroth-order schemes based only on function-value feedback, which further increases the complexity of the problem.}

Motivated by these considerations, this paper studies resilient distributed online optimization over time-varying directed networks with unknown Byzantine edge attacks. Moreover, we consider the more challenging bandit feedback setting, aiming to achieve consensus and sublinear regret in dynamic adversarial environments.
\begin{table*}[t]
\centering
%\color{blue}
\caption{Comparison of Distributed Optimization Methods under Byzantine Attacks}
\label{tab:comparison}
\begingroup
\fontsize{7pt}{8.4pt}\selectfont
\setlength{\tabcolsep}{2.6pt}
\renewcommand{\arraystretch}{1.13}
\begin{tabularx}{\textwidth}{@{}ccccc>{\raggedright\arraybackslash}X>{\raggedright\arraybackslash}Xc@{}}
\toprule
Reference & Setting & \makecell{Topology,\\Attack Set$^{a}$}
& Target & \makecell{Zeroth-\\order} & Optimality/Regret
& {Key Conditions$^{c}$}
& \makecell{Exact Optimum /\\Sublinear Regret$^{d}$} \\
\midrule
\makecell[tc]{\cite{2019Sundaram}/\cite{2020Su}\\/\cite{2019Zhao} }& Offline & Static,Fixed & Node & No
& %Approach to 
Local-minimizer hull
& Graph robustness%; convexity 
/Reduced-graph/Connected dominating trusted agents set
& No \\
\addlinespace[2pt]
\cite{2022Zhu} & Offline & Static,Fixed & Node & No
& %Exact convergence to the 
Original optimum
& Graph and loss redundancy %nonempty-interior optimal set 
& Yes \\
\addlinespace[2pt]
\cite{2025Kuwaranancharoen} & Offline & Static,Fixed & Node & No
& %Geometric approach to an 
Optimum neighborhood%; approximate consensus
& %Strong convexity; %smoothness; 
Contraction,mixing conditions 
& No \\
\addlinespace[2pt]
\cite{2023Xu2}/\cite{2023Xu} & Offline & \makecell{Static, fixed\\/Local bound$^{b}$}
& \makecell{Node \&\\Edge} & No
& %Consensus; 
local-minimizer interval
& %Per-round local attack bound; 
Prescribed trusted-source/Graph robustness
& No \\
\midrule
\cite{2022Sahoo} & Online & Static,Fixed & Node & No
& %$O((\log T)^2)$ 
Sublinear weighted static regret
& %Local strong convexity; 
Graph robustness & No \\
\addlinespace[2pt]
\cite{2024Wei} & Online & Static,Fixed & Node & Yes
& Sublinear expected weighted dynamic regret
& Graph robustness; path 
variation% and exploration conditions 
& No \\
\addlinespace[2pt]
\cite{2023Dong} & Online & Static,Fixed & Node & No
& Linear adversarial static regret; sublinear stochastic regret
& Central-server setting; i.i.d. losses and regularity conditions %for stochastic sublinearity 
& No \\
\midrule
\textbf{This work} & \textbf{Online}
& \makecell[tc]{\textbf{Time-varying,}\\\textbf{Time-varying}}
& \textbf{Edge} & \textbf{Yes}
& \textbf{Sublinear standard dynamic regret for the loss sum}
& \makecell[tl]{Trusted links with protected scalar \\feedback;
path variation}% and parameter choices
& \textbf{Yes} \\
\bottomrule
\end{tabularx}

\vspace{4pt}
\begin{minipage}{\textwidth}
\fontsize{7pt}{8.4pt}\selectfont
\textit{Notes.} 
\textit{(a)} Topology and attack set refer to the communication network and the identities of adversarial nodes or compromised edges, respectively. 
\textit{(b)} ``Local bound'' means the number of altered incoming messages at each iteration is bounded in \cite{2023Xu};  
\textit{(c)} Only distinguishing conditions are listed; %the full assumptions are given in the cited works. 
\textit{(d)} For offline methods, ``Yes'' indicates exact convergence to a minimizer of the original unweighted objective; For online methods, it indicates sublinear standard static or dynamic regret for the original unweighted loss sum. 
%The unweighted benchmark sums nonfaulty agents' costs for node attacks and all agents' costs for edge attacks. Expectations over algorithmic randomization are permitted. Population-risk stochastic regret, as in~\cite{2023Dong}, uses a different benchmark and is not counted as standard empirical regret here. 
``No'' means this specified guarantee is not established. 
%, not that it is impossible in every special case. For this work, sublinearity requires the stated parameter choices and $\sum_{t=1}^{T}\|x_{t+1}^{*}-x_t^{*}\|/\alpha_t=o(T)$.
\end{minipage}
\endgroup
\end{table*}

\subsection{Literature review}

Distributed convex optimization with time-invariant objective functions has been extensively studied in recent years. In this framework, each agent has access only to its own local convex objective and cooperates with its neighbors to minimize a global objective, typically the sum of local functions. Starting from the seminal work~\cite{2009Nedic}, a large body of results has been developed; see, e.g.,~\cite{Lobel2010,2014Chang,2010Ram} and references therein. Zeroth-order distributed optimization has also been investigated as an alternative when gradient information is unavailable~\cite{2019Wang,2020Tang}. However, these works typically assume reliable communication and do not account for adversarial attacks. 

%Resilient distributed optimization under adversarial attacks has received growing attention. {For Byzantine \emph{node} attacks, local filtering–based strategies have been proposed to mitigate adversarial information~\cite{2019Sundaram,2020Su}, %These methods establish asymptotic agreement together with approach to the convex hull of the nonfaulty agents' local minimizers~\cite{2019Sundaram}, or to the generally smaller union of minimizer sets of admissible weighted combinations of their local costs~\cite{2020Su}.
%both of which guarantee the nonfaulty agents' estimates asymptotically approach the convex hull of their local minimizer sets. Trusted-agent-based methods provide convergence to an admissible weighted-optimal set associated with trusted agents~\cite{2019Zhao}. For strongly convex and smooth costs, geometric convergence to a neighborhood of the nonfaulty agents’ equal-weight aggregate minimizer, together with approximate consensus, is established in~\cite{2025Kuwaranancharoen}. }{Note that above guarantees concern convex-hull confinement, weighted optimality, or neighborhood convergence, rather than general exact minimization of the unweighted sum of the nonfaulty agents' local costs. }{By contrast, exact convergence is established in~\cite{2022Zhu} under restrictive structural graph and objective redundancy assumptions: the sum of costs over every sufficiently large agent subset must have the same minimizer set as the full aggregate. }

Resilient distributed optimization under adversarial attacks has
received growing attention. {Under Byzantine \emph{node} attacks, local-filtering methods drive
nonfaulty agents' estimates toward the convex hull of their local minimizer sets~\cite{2019Sundaram,2020Su}, while trusted-agent methods converge to a minimizer of an admissible weighted sum of the trusted
agents' objectives~\cite{2019Zhao}. Under strong convexity and smoothness conditions, \cite{2025Kuwaranancharoen} establishes approximate
consensus and geometric convergence to a neighborhood of a minimizer of the sum of the nonfaulty agents' objectives. Zhu et al.~\cite{2022Zhu} show that all nonfaulty agents converge to a common minimizer of that sum. Their result requires the communication graph to remain rooted after specified node and edge removals and the set of minimizers of the summed objectives to remain unchanged when any subset of agents meeting a graph-dependent size threshold is retained.}

In comparison, Byzantine \emph{edge} attacks have received much less attention. {Distance-based and trusted-state interval filtering are proposed in~\cite{2023Xu} and~\cite{2023Xu2}, respectively, with guarantees of agreement and confinement to intervals determined by local minimizers, rather than general exact
minimization of the original aggregate objective. Under weaker, structured models of bounded false-data injection, \cite{2023Gusrialdi,2023Iqbal} establish neighborhood convergence. %{Under more structured false-data-injection models, \cite{2023Gusrialdi} considers bounded attacks on both local computations and communication links, while \cite{2023Iqbal} studies bounded, structured communication injections through competitive-interaction dynamics.}  
These methods establish neighborhood optimality under their respective assumptions. However, all the aforementioned studies are restricted to time-invariant local objective functions.}

To address time-varying environments, distributed online convex optimization (DOCO) has been introduced as an extension of online convex optimization~\cite{2016Hazan}. In DOCO, agents update decisions based on revealed losses without access to future objectives and  sublinear regret is a standard performance criterion. While sublinear static regret is widely obtained~\cite{2013Yan,2013Hosseini,MateosNunez2014}, with zeroth-order extensions studied in~\cite{2021Yi,2023Liu,2022Xiong}, it is insufficient to capture tracking performance in dynamic environments~\cite{2007Hazan}. To address this limitation, dynamic regret~\cite{2018Shahrampour,2019Zhang,2019Lu,2022Xiong} has been proposed as a stronger metric. Although more challenging to analyze, it provides a more appropriate characterization of performance and is therefore adopted in this paper.

Only a few works study distributed {online} optimization under Byzantine attacks, and these works mainly address \emph{node} attacks. {%In decentralized networks,
Sahoo et al.~\cite{2022Sahoo} establish $O((\log T)^2)$ weighted static regret under locally strongly convex losses, whereas \cite{2024Wei} obtains sublinear expected weighted dynamic regret under path-variation and exploration-parameter conditions. Both regret criteria involve mixing-induced, generally unavailable node weights, replacing the loss sum with a convex combination; our method avoids this relaxation. In a central-server setting, \cite{2023Dong} considers honest agents' unweighted loss sum, establishing a tight linear static-regret bound and, under i.i.d. objectives, sublinear stochastic regret against a fixed population-risk minimizer. }%{, which is a relatively strong assumption for online problems.}
%\cite{2023Dong} analyzes regret with respect to the equal-weight average of honest participants’ losses: its robust online gradient method admits a linear adversarial static-regret term that is attained by constructed examples, while a momentum method achieves sublinear stochastic regret under independent losses from a common distribution. The stochastic and adversarial guarantees use different benchmarks

To the best of our knowledge, distributed online optimization under Byzantine \emph{edge} attacks has not been studied before. Moreover, no existing work establishes sublinear standard regret  %{ defined using the unweighted sum of all agents' local losses} 
without imposing stochastic assumptions (e.g., i.i.d.) on the loss functions under Byzantine attacks. %{Our analysis does not require i.i.d. losses and establishes variation-dependent  sublinear dynamic regret bounds.} 
Table~1 summarizes the comparison with related work. 

%Compared with the node-attack literature, distributed online optimization under Byzantine \emph{edge} attacks remains less explored. We study this setting with regret defined using the unweighted sum of all agents' local losses. Our analysis does not require i.i.d. losses and establishes variation-dependent dynamic regret bounds that are sublinear under explicit conditions on the minimizer path and algorithm parameters. Table~1 summarizes the comparison with related work.

\subsection{Contributions}
The main contributions of this paper are as follows. 
\begin{itemize}
\item[a)] %We formulate a distributed online optimization problem under time-varying Byzantine edge attacks, {allowing both the communication topology and the set of compromised edges to vary over time.} This setting captures more realistic and dynamic security threats than existing models focusing on fixed adversarial-node sets or time-invariant objectives.
We formulate distributed online optimization under Byzantine
edge attacks with time-varying communication topology and
compromised-edge sets, capturing dynamic threats beyond models
with fixed adversarial-node sets or time-invariant objectives. 

\item[b)] We develop two classes of zeroth-order distributed online optimization algorithms resilient to Byzantine attacks, based on randomized and deterministic difference, %schemes, 
respectively. %, making them applicable when explicit gradients are unavailable. 
In particular, numerical results indicate the deterministic difference scheme leads to faster consensus formation, enhanced post-consensus stability, and improved convergence performance.
%In the reported simulations, the deterministic scheme exhibits faster consensus and smaller post-consensus fluctuations.

\item[c)] %We establish sublinear dynamic regret bounds for the proposed algorithms. {Unlike the weighted-regret analyses and the i.i.d.-based stochastic-regret result for node-attack models, our guarantees do not introduce unknown node weights or require i.i.d. losses.% and do not require the objective-redundancy assumptions used in~\cite{2022Zhu}. 
We establish sublinear dynamic regret under the stated path-variation and parameter conditions. {Unlike weighted-regret or i.i.d.-based stochastic-regret results under node attacks, our guarantees do not introduce unknown node weights or require i.i.d. losses. The bounds exhibit horizon and path-variation dependence comparable to attack-free zeroth-order guarantees under corresponding parameter choices~\cite{2022Xiong}.}
\end{itemize}

The rest of this paper is organized as follows. In Section~2, the preliminaries, problem formulation and the algorithms are presented. In Section~3, consensus of the agents' estimates and regret bounds are established. Simulation examples are provided in Section~4, and conclusions are drawn in Section~5. 

{\em Notation.} For a matrix $M$, let $M[i,j]$, $M[i,:]$, and $M[:,j]$ denote its $(i,j)$-th entry, $i$-th row, and $j$-th column, respectively. Let $\textbf{e}^k\in\mathbb{R}^d$ be the $k$-th canonical basis vector. Denote by $\|\cdot\|$ and $\langle\cdot,\cdot\rangle$ the Euclidean norm and inner product on $\mathbb{R}^d$, respectively. The projection onto $X$ is denoted by $\mathcal{P}_X(y)\triangleq \arg\min_{x\in X}\|x-y\|$. A vector is called stochastic if its entries are nonnegative and sum to one; a matrix is row-stochastic if each row is stochastic, and doubly stochastic if both $M$ and $M^\top$ are row-stochastic. Let $\mathbb{B}$ and $\mathbb{S}$ denote the unit ball and unit sphere in $\mathbb{R}^d$, respectively. %{Denote $\mathrm{Unif}(W)$ the uniform distribution over a set $W$. }

\section{Preliminaries, Problem Formulation and Algorithm Design}

\subsection{Preliminaries: zeroth-order DOCO}
We consider a multi-agent system with agents indexed by $V=\{1,\dots,n\}$. The agents cooperatively aim to minimize a time-aggregated global loss $\sum_{t=1}^T f_t(x_t)$ over a $T$-horizon, where at each time $t\ge 0$,
\begin{align*}
    f_t(x)=\sum_{i=1}^n f_t^i(x),
\end{align*}
where $f_t^i:X\to\mathbb{R}$ is the local loss known only to agent $i$, and $X\subset\mathbb{R}^d$ is a common feasible set.
\paragraph*{Online and bandit feedback} 

The local losses may vary arbitrarily over time. At time $t$, agent $i$ has access to $f_1^i,\dots,f_t^i$ but NOT to the future loss $f_{t+1}^i$ when choosing $x_{t+1}^i$. Moreover, we consider \emph{bandit-feedback} setting: $i$ can query values of $f_t^i(x)$ at selected points $x\in X$, but cannot directly access gradient information.
\paragraph*{Performance metrics} 
A standard performance metric is the \emph{static regret} of agent $i$~\cite{MateosNunez2014},
\begin{align*}
\mathrm{regret}_{T}^{S,i}
= \sum_{t=1}^{T}\Big(f_t(x_t^i)-f_t(x^*)\Big),
\end{align*}
where $x^*\in\arg\min_{x\in X}\sum_{t=1}^T f_t(x)$ is the best fixed decision in hindsight. Since static regret may fail to characterize tracking performance in nonstationary environments~\cite{2007Hazan}, we adopt the \emph{dynamic regret}~\cite{2019Lu,2022Xiong}:
\begin{equation}\label{5}
\mathrm{regret}_{T}^{D,i}
= \sum_{t=1}^{T}\Big(f_t(x_t^i)-f_t(x_t^*)\Big),
\end{equation}
where $x_t^*\in\arg\min_{x\in X} f_t(x)$ is a per-round minimizer.
%{\begin{remark}
%Unlike static optimization with a fixed objective, online optimization requires agents to make decisions without knowing future loss functions and to adapt to a potentially time-varying global minimizer. Accordingly, dynamic regret is used to quantify the cumulative tracking error, and its analysis must account for the temporal variation of the minimizer sequence. In our setting, this challenge is compounded by time-varying communication networks, Byzantine edge attacks, and zeroth-order feedback. When the objective is time-invariant, the minimizer-path variation term vanishes, and the analysis reduces to the corresponding static-objective case.
%\end{remark}}

\paragraph*{Communication graph}

Agents communicate over a time-varying directed graph $\mathcal{G}_t=(V,E_t)$, where $(i,j)\in E_t$ means that agent $j$ sends its current decision to agent $i$ at time $t$. Let $\mathcal{N}_t^i=\{j\in V:(i,j)\in E_t\}$ denote the in-neighbor set of agent $i$. A directed graph is strongly connected if there exists a directed path between every ordered pair of nodes. 
The details of the information exchange protocol will be specified later.

\subsection{Problem Formulation: Zeroth-order DOCO under Byzantine edge attacks}

We consider Byzantine edge attacks on the time-varying directed graph $\mathcal{G}_t=(V,E_t)$. At each time $t\in\{1,\dots,T\}$, the edge set is partitioned as
\begin{align*}
    E_t = E_t^{r} \cup E_t^{n} \cup E_t^{a},
\end{align*}
where $E_t^{r}$, $E_t^{n}$, and $E_t^{a}$ denote the trusted, normal, and adversarial edges, respectively. Messages transmitted over adversarial edges may be arbitrarily manipulated, whereas those over trusted and normal edges remain truthful. {We assume that each agent knows which of its incoming links are trusted, whereas its other incoming links cannot be distinguished
as normal or adversarial.} In addition, each agent has a trusted self-loop, i.e., $(i,i)\in E_t^{r}$ for all $i\in V$. {More precisely, the edge-dependent message received by agent $i$ from sender $j$ is}
\begin{align}
{\tilde x_t^{ij}=\begin{cases}
x_t^j,&\text{if }\,(i,j)\in E_t^r\cup E_t^n,\\
a_t^{ij},&\text{if }\,(i,j)\in E_t^a,
\end{cases}}
\label{eq:recv}
\end{align}
{where $a_t^{ij}\in\mathbb R^d$ is arbitrary and receiver dependent. }

We impose the following standard assumptions. 
\begin{assumption}\label{asmp:X-interior}
The feasible set $X\subset\mathbb{R}^d$ is compact and convex, and contains a ball of radius $r>0$, i.e., $r\mathbb{B}\subset X$.
\end{assumption}

Assumption~\ref{asmp:X-interior} is standard in zeroth-order optimization~\cite{2005Flaxman} ensuring the feasibility of random perturbations. {Define the finite radius $M\triangleq\max_{x\in X}\|x\|$.}
\begin{assumption}\label{asmp:convex-feasible-set}
For any $t\geqslant 0$ and $i\in V$, $f_t^i$ is convex and differentiable on $X$, and there exists $L>0$, independent of $i$ and $t$, such that $\|\nabla f_t^i(x)\|\le L$ for any $x\in X$.
\end{assumption}

We impose the following assumption to ensure sufficient exchange of trusted information.

\begin{assumption}\label{asmp:strongly-connected-over-time}
There exists a positive integer $ B_0$ such that the graph $\left(V,E_t^{r} \cup E_{t+1}^{r}\cup \ldots \cup E_{t+B_0-1}^{r}\right)$ is strongly connected for any $ t \geqslant 0$.  
{Moreover, for every $(i,j)\in E_t^r$, a protected channel logically isolated from state-message transmission reliably delivers uncorrupted scalar feedback from $i$ to $j$ within round $t$.}
\end{assumption}

%{Under these assumptions, our objective is to design zeroth-order distributed algorithms that (i) reaches consensus, and (ii) controls the dynamic regret in~\eqref{5} for the original unweighted aggregate loss $f_t=\sum_{i=1}^n f_t^i$. Every agent and every local loss remains legitimate; only selected transmissions are compromised.}
{
\begin{remark}[Trusted communication] 
%Trusted communication can be supported by physically isolated, access-controlled links, with their deployment and scheduling ensuring joint connectivity. Protected feedback can be supported by a protected out-of-band management 
%and control-message protection 
%mechanism~\cite{RFC4778}. %,RFC9055}. If the trusted graph is undirected, i.e., $(i,j)\in E_t^r$ implies $(j,i)\in E_t^r$, feedback can be transmitted over the reverse trusted edges, and no separate feedback channel is needed.
Assumption~3 models a managed network with selectively protected
communication links. Trusted communication can be supported by physically isolated, access-controlled links, whose deployment
and scheduling are designed to ensure joint connectivity. Scalar feedback can be carried over a protected out-of-band control channel~\cite{RFC4778}. If the trusted graph is undirected, i.e., $(i,j)\in E_t^r\iff (j,i)\in E_t^r$, the reverse trusted edges can carry the protected scalar feedback without additional physical feedback links.
\end{remark}
}

\subsection{Algorithm Design}
%In this paper, we propose two DOCO algorithms, both of which consist of three main mechanisms: 
%1) neighbor filtering; 
%2) construction of attack-resilient mixing coefficients and update rule;
%3) zeroth-order gradient estimation. 
%The two algorithms differ only in the third mechanism.

In this paper, we propose two DOCO algorithms, both of which share the same neighbor filtering and resilient mixing--descent steps, and differ only in the construction of the zeroth-order direction.

\subsubsection{Neighbor filtering} 

At time $t\in\{1,\dots,T\}$, for each agent $i$, $j$ sends $x_t^j$ to $i$ over each edge $(i,j)\in E_t$, {but $i$ receives $\tilde{x}_t^{ij}$ defined in~\eqref{eq:recv} since compromised incoming links are unknown. Each agent first filters received messages. }

Inspired by~\cite{2019Zhao}, for $i\in V,~t\ge 0$, we adopt a coordinate-wise neighbor filtering scheme. At time $t$, agent $i$ constructs a trusted range using only the trusted incoming edges $E_t^r$. For each coordinate $k\in\{1,\dots,d\}$, define
\begin{align*}
  M_t^i[k] = \max_{(i,j)\in E_t^r}\{x_t^j[k]\}, 
  m_t^i[k] = \min_{(i,j)\in E_t^r}\{x_t^j[k]\}.
\end{align*}
Agent $i$ then forms the resilient neighbor set
\begin{align}
S_t^i
\!=\!\Bigl\{j\!\in\! \mathcal N_t^{i}\!\cup\!\{i\}:
\tilde x_t^{ij}\!\in\! \prod_{k=1}^d[m_t^i[k],M_t^i[k]]{\cap X_{\delta_{t\!-\!1}}}\Bigr\}\label{eq:S_t^i}
\end{align}
{for $t\ge 0$, where $\delta_{-1}\triangleq \delta_0$. }  
Neighbors in $S_t^i$ are referred to as \emph{safe neighbors}, while those in $(\mathcal N_t^{i}\cup\{i\})\setminus S_t^i$ are regarded as \emph{unsafe neighbors} and are excluded from agent $i$'s update at time $t$.
% Define the \emph{resilient neighbor set} by

%\begin{align}
%{S_t^i
%=\Bigl\{j\in \mathcal N_t^{i}\cup\{i\}:
%\tilde{x}^{ij}_t\in \mathcal C_t^i\Bigr\}}.  \label{eq:S_t^i}
%\end{align}
%Neighbors in $S_t^i$ are referred to as \emph{safe neighbors}, while those in $\mathcal N_t^{i}\setminus S_t^i$ are regarded as \emph{unsafe neighbors} and excluded from $i$'s update at time $t$. 

\subsubsection{Attack-resilient mixing and update rule}
Both algorithms update $x_t^i$ as
\begin{equation}\label{update}
x_{t+1}^i=\mathcal{P}_{X_{\delta_t}}\left(
y_t^i-\alpha_t\frac{G_t^i}{z_{t+1}[i]}
\right),
\end{equation}
where $\alpha_t>0$, $\delta_t\in(0,r]$, and $X_{\delta_t}\triangleq(1-\delta_t/r)X$. Here $y_t^i$ is the resiliently mixed state, $z_t[i]$ is a push-sum scaling variable, and $G_t^i$ is a zeroth-order direction deferred to the next subsection. Set
\begin{align}
    y_t^i \;=\; \sum_{j=1}^n\tilde{A}_t[i,j]x_t^j,\label{y}
\end{align}
where $\tilde{A}_t[i,j]$ is the \emph{attack-resilient mixing coefficient}, {containing the reconstructed messages of all safe neighbors in $S_t^i$. 
Inspired by the push-sum methods in~\cite{2014Nedic,2019Scutari,2022Xiong},
we design a mixing rule that uses trusted links to mitigate Byzantine edge attacks. To construct $\tilde{A}_t[i,j]$, we first assign a weight
$s_t[i,q]$ to each trusted link $(i,q)\in E_t^r$. We then normalize these weights across the trusted outgoing links of each sender $q$ to obtain $B_t[i,q]$.

By \eqref{eq:S_t^i}, neighbors connected by trusted links belong to $S_t^i$. Next, we reconstruct the accepted nontrusted messages as convex combinations of trusted states, ensuring that subsequent mixing uses only actual agent states. Define the \emph{trusted convex-hull} by
\begin{align}
\begin{aligned}[b]
&\mathcal C_t^i
\triangleq \operatorname{co}\{x_t^q:(i,q)\in E_t^r\}=\Biggl\{
\sum_{q:(i,q)\in E_t^r}\lambda_q x_t^q:
\\&
\lambda_q\!\geq\!0\ \text{for all }q\text{ with }(i,q)\!\in\! E_t^r,
\sum_{q:(i,q)\in E_t^r}\!\lambda_q\!=\!1\!
\Biggr\}, 
\end{aligned}
\label{eq:filtering}
\end{align}
where  $\operatorname{co}(\cdot)$ denotes convex hull. Since $(i,i)\in E_t^r$, $\mathcal C_t^i\neq \emptyset$. 
Define the coefficient $\lambda_t^{ij}
%\in \mathbb{R}^{|{q:(i,q)\in E_t^r}|}
$ for $j\in S_t^i,~(i,j)\notin E_t^r$ by
\begin{align}
\nonumber&\lambda_t^{ij}=\operatorname{argmin}_{\lambda\in \Lambda_t^{ij}}\|\lambda\|^2,\text{where }
\Lambda_t^{ij}=\Big\{\lambda\!\in \!\mathbb{R}^{|{q:(i,q)\in E_t^r}|}
\\&\!:\!\lambda[q]\!\geq\!0, 
\!\sum_{q:(i,q)\in E_t^r}\!\!\lambda[q]\!=\!1,\mathcal{P}_{\mathcal{C}_t^i}(\tilde{x}^{ij}_t)\!\!=\!\!\!\sum_{q:(i,q)\in E_t^r}\!\!\!\lambda[q] x_t^q\!
%\in \!\mathcal{C}_t^i
\Big\}\label{eq:lambda}
\end{align}
Since $\Lambda_t^{ij}$ is nonempty, compact, and convex,
strict convexity of $\|\lambda\|^2$ ensures the minimizer in \eqref{eq:lambda} exists and is unique. 
Both the projection and coefficient selection in \eqref{eq:filtering}--\eqref{eq:lambda} reduce to linearly constrained convex quadratic programs, which are polynomial-time solvable ~\cite{1980Kozlov}.%under the standard rational-input bit-complexity model~\cite{1980Kozlov}.}
\begin{remark}[Selection of coefficient $\lambda_t^{ij}$] % \eqref{eq:filtering}--\eqref{eq:lambda} gives one setting of the coefficient $\lambda_t^{ij}$. The way of presenting the coefficients is not limited to this method. It needs to be selected based on the dimension, the distribution of trusted edges, and computing resources, etc. As long as $\lambda[q]\geq 0,\sum_{q:(i,q)\in E_t^r}\lambda[q]=1$ is guaranteed, it will not affect the subsequent theoretical analysis.
\eqref{eq:filtering}--\eqref{eq:lambda} provide one construction of $\lambda_t^{ij}$. Alternative local constructions may be chosen according to the problem dimension, trusted-link structure, and computational resources. The main theoretical conclusions remain unchanged provided
$\lambda_t^{ij}[q]\geq0$ for all neighbors $q$ linked by trusted edges and $\sum_{q:(i,q)\in E_t^r}\lambda_t^{ij}[q]=1$.
\end{remark}
We define the weights $s_t[i,q]$ as follows:
\begin{equation}
s_t[i,q]\triangleq
\begin{cases}
1+\displaystyle\sum_{j:~j\in S_t^i,(i,j)\notin E_t^r}
\lambda_t^{ij}[q],&(i,q)\in E_t^r,\\[1ex]
0,&(i,q)\notin E_t^r.
\end{cases}
\label{eq:s[i,q]}
\end{equation}
The constant $1$ is the baseline contribution of trusted message delivered by $q$ itself. Each accepted nontrusted message contributes one additional unit, distributed among the trusted anchors by its barycentric coefficients. Thus $s_t[i,q]$ is receiver $i$'s score for anchor $q$. The exact absorption identity is
\begin{equation*}
\sum_{q:(i,q)\in E_t^r}x_t^q
\!+\!\sum_{j:~j\in S_t^i,(i,j)\notin E_t^r}\!\mathcal P_{\mathcal C_t^i}(\tilde x_t^{ij})\!=\!\sum_{q=1}^n s_t[i,q]x_t^q.
\end{equation*}
Since at most $n-1$ nontrusted messages can be accepted, 
\begin{equation}
1\leq s_t[i,q]\leq n
\quad\text{whenever }(i,q)\in E_t^r.
\label{score-bounds}
\end{equation}
Note that $\mathcal C_t^i$, $S_t^i$,  $\lambda_t^{ij}$ and $s_t[i,q]$ can be locally calculated by node $i$. 

Then we can give the column normalization $B_t[i,q]$. By Assumption \ref{asmp:strongly-connected-over-time}, scalar $s_t[i,q]$ can be safely delivered to agent $q$ by $i$ for all $i:~(i,q)\in E_t^r$. Node $q$ set $s_t[j,q]=0$ for $j$ that no $s_t[j,q]$ delivered in this turn. Then all agent $q\in V$ can calculate 
%q接收到所有i的信息
\begin{equation}
%\kappa_t[q]\!\!\triangleq\!\!\!\!\sum_{i:(i,q)\in E_t^r}\!\!\!s_t[i,q]\!=\!\sum_{i=1}^n\! s_t[i,q],~B_t[i,q]\!\!\triangleq\!\!\frac{s_t[i,q]}{\kappa_t[q]}.
\kappa_t[q]\triangleq\sum_{i=1}^n s_t[i,q],~B_t[i,q]\triangleq\frac{s_t[i,q]}{\kappa_t[q]}.
\label{B-construction}
\end{equation}
The trusted self-loop gives $\kappa_t[q]\geq 1$, so the normalization is well defined. Moreover, we have 
\begin{equation}
\sum_{i=1}^nB_t[i,q]=1,
~~
B_t[i,q]>0\Longleftrightarrow(i,q)\in E_t^r.
\label{eq:B-column}
\end{equation}
\eqref{score-bounds}-\eqref{B-construction} yield $\kappa_t[q]\leq n^2$, and each trusted entry satisfies 
\begin{equation}
B_t[i,q]\geq n^{-2}\triangleq\beta_B,~(i,q)\in E_t^r
\label{B-lower}
\end{equation}
Initialize $z_0[i]=1/n$ for all $i\in V$. At time $t\ge 0$, each agent $q$ transmits $B_t[i,q]z_t[q]$ to $i$ through $(i,q)\in E_t^r$. For $q:~(i,q)\notin E_t^r$, $i$ sets $B_t[i,q]z_t[q]=0$. Then, $i$ updates the \emph{push-sum scaling variable} by 
\begin{align}
z_{t+1}[i]=\sum_{q=1}^nB_t[i,q]z_t[q], \label{z}
\end{align}
Since $B_t$ is column-stochastic, $z_t$ remains stochastic. Accordingly, define $\tilde{A}_t\in\mathbb R^{n\times n}$ by 
\begin{equation}\label{A}
\tilde A_t[i,j]\triangleq
\frac{B_t[i,j]z_t[j]}{z_{t+1}[i]}. 
\end{equation}
Equations~\eqref{z}--\eqref{A} give
\begin{align}
\sum_{j=1}^n\tilde A_t[i,j]&=1,
\label{A-row}\\
z_{t+1}^{\top}\tilde A_t&=z_t^{\top}.
\label{row}
\end{align}
%[Role of accepted nontrusted messages]
%The proposed filtering procedure does not necessarily identify all corrupted messages. Even if an adversarially manipulated message \(\widetilde{x}_t^{ij}\neq x_t^j\) belongs to \(S_t^i\), it is not directly used in the state-mixing update. Instead, its projection onto the trusted convex hull \(C_t^i\) is employed only to construct barycentric coefficients, which are subsequently incorporated into the trusted-based weights \(s_t[i,q]\). As a result, the mixing matrix \(\widetilde A_t\) has nonzero entries only on trusted communication edges, and the mixed state \(y_t^i\) is a convex combination of true states received from trusted neighbors. Hence, Byzantine edge attacks may influence the mixing weights, but cannot directly inject corrupted state vectors into the resulting recursion.
%The filtering procedure need not identify all corrupted messages. An accepted nontrusted message is projected onto the trusted convex hull and used only to construct barycentric coefficients for the trusted-edge weights $s_t[i,q]$. Thus, $\widetilde A_t$ mixes only true states: attacks may affect the mixing weights but cannot directly inject corrupted states into the recursion. Moreover, the numerical results in Section 4 provide empirical evidence that exploiting such messages improves convergence speed compared with using trusted edges alone. 
\begin{remark}[Edge attacks and resilient mixing]\label{rem:role}
Unlike Byzantine node attacks, which can compromise local computations, Byzantine edge attacks in our model corrupt messages in transit. Every agent evaluates its own local loss and follows the prescribed updates, so the global objective includes all agents’ losses with equal weight. By projecting accepted messages from nontrusted links onto the convex hull of states received over trusted links, we express each mixed state as a convex combination of actual agent states. Corrupted messages can thus influence the mixing weights but enter the update only through those weights. The push-sum correction cancels the node-dependent weights in the aggregate descent bound, yielding a dynamic regret guarantee for the original unweighted objective.
\end{remark}
}
\begin{framework}[t]
\caption{A Unified Framework for Zeroth-Order Resilient DOCO under Byzantine Edge Attacks}
\label{fr:common}
\begin{algorithmic}[1]
\Require Horizon $T$, stepsizes $\{\alpha_t\}_{t=0}^T$, perturbation radii $\{\delta_t\}_{t=0}^T$, and initial estimates $\{x_0^i\}_{i\in V}\subset X_{\delta_0}$, $z_0=[1/n,\ldots,1/n]^\top$. 
\For{$t=0,\ldots,T$ %and $i=1,\ldots,n$
}
\State {Each $i$ receives $\tilde{x}_t^{ij}$ in \eqref{eq:recv} from neighbors, computes the \textbf{resilient neighbor set $S_t^i$} by \eqref{eq:S_t^i}.} %{and the \textbf{%filtering  trusted convex-hull} $\mathcal{C}_t^i$ by \eqref{eq:filtering}.}
\State {Each $i$ computes %the \textbf{trusted-based safe weight }
$s_t[i,q]$ for $\{q:~(i,q)\in E_t^r\}$ by \eqref{eq:filtering}-\eqref{eq:s[i,q]} and sends each score to $q$ through the protected scalar-feedback channel.}
\State {Each node $q%~(\exists i,~s.t.(i,q)\in E_t^r)
$ collects $s_t[i,q]$, computes $\kappa_t[q],~B_t[i,q]$ by~\eqref{B-construction} and transmits $B_t[i,q]z_t[q]$ to $i$. 
}
\State {Each $i$ receives $B_t[i,q]z_t[q]$ from trusted $q$ and sets it to be $0$ otherwise, updates \textbf{push-sum scaling variable $z_{t+1}[i]$} by \eqref{z}, and computes the \textbf{resiliently mixed state $y_t^i$} by \eqref{y} and \eqref{A}}
\begin{align}
{y_t^i=\sum_{q=1}^n\frac{B_t[i,q]z_t[q]}{z_{t+1}[i]}x_t^q.}\label{sa1}
\end{align}
\State Each $i$ constructs the \textbf{zeroth-order direction $G_t^i$} as specified in Algorithms~\ref{alg:rand}--\ref{alg:det} and updates \textbf{$x_{t+1}^i$} by \textbf{difference-based descent:}
\begin{align}
x_{t+1}^i=\mathcal{P}_{X_{\delta_t}}\left(y_t^i-\alpha_{t}\frac{G_t^i}{z_{t+1}[i]}\right).\label{sa4}
\end{align}
\EndFor
\end{algorithmic}
\end{framework}
\subsubsection{Zeroth-order gradient estimation}
Next, we specify $G_t^i$ using two zeroth-order difference schemes, randomized and deterministic. 
\paragraph*{Randomized difference via $\delta$-smoothing}
Following~\cite{2005Flaxman}, define
\begin{align}
\widehat{f}^i_{t,\delta}(x)\triangleq\mathbb{E}[f^i_t(x+\delta v)],\label{aa1}
\end{align}
where $\delta\in(0,r]$ and $v$ is uniformly distributed over $\mathbb{B}$. It
can be directly proved that $X_{\delta}+\delta\mathbb{B}\subset X$, 
and thus $\widehat{f}^i_{t,\delta}(\cdot)$ is well-defined on $X_{\delta}$.  
Moreover, it holds that (see~\cite{2005Flaxman})
\begin{align}
\nabla\widehat{f}^i_{t,\delta}(x)=\frac{d}{\delta}\mathbb{E}[f^i_t(x+\delta u)u], \label{aa1'}
\end{align}
{where $u$ is uniformly distributed over $\mathbb{S}$, denoted by $u_t^i\sim \mathrm{Unif}(\mathbb S)$}, and equivalently
\begin{align}
\nabla\widehat{f}^i_{t,\delta}(x)
=\frac{d}{2\delta}\mathbb{E}\!\left[\big(f^i_t(x+\delta u)-f^i_t(x-\delta u)\big)u\right].\label{aa2'}
\end{align}
Hence, drawing $u_t^i\sim \mathrm{Unif}(\mathbb S)$, the sample 
\begin{align*}
g_t^i=\frac{d}{2\delta_t}\left[f^i_t(y_t^i+\delta_t u_t^i)-f^i_t(y_t^i-\delta_t u_t^i)\right]u_t^i
\end{align*}
is an unbiased estimator of $\nabla\widehat{f}^i_{t,\delta_t}(y_t^i)$. 
Therefore, in~\eqref{update} we set $G_t^i=g_t^i$, which yields Algorithm~\ref{alg:rand}.
\paragraph*{Deterministic difference}
%To improve convergence speed, we also consider a deterministic perturbation scheme motivated by the Kiefer--Wolfowitz method~\cite{1952Kiefer}. Assume that $f_t^i(\cdot)$ can be queried at $y_t^i\pm \delta_t \textbf{e}^k$, $k=1,\ldots,d$, $\delta_t>0$, and denote $[w_t^{i,k}]^{+}=f_t^i(y_t^i+\delta_t\textbf{e}^k)$ and $[w_t^{i,k}]^{-}=f_t^i(y_t^i-\delta_t\textbf{e}^k)$. Define 
Motivated by the Kiefer--Wolfowitz method~\cite{1952Kiefer}, query $f_t^i$ at $y_t^i\pm \delta_t \textbf{e}^k$, obtain $[w_t^{i,k}]^{\pm}$ for $k=1,\ldots,d$, and define
\begin{align*}
&h_t^i[k]=\frac{[w_t^{i,k}]^{+}-[w_t^{i,k}]^{-}}{2\delta_t},~k=1,\cdots,d,\\
&h_t^i = \left[h_t^i[1],\cdots,h_t^i[d]\right]^{\top},
\end{align*}
which serves as the deterministic difference estimator. 
Accordingly, in~\eqref{update} we set $G_t^i=h_t^i$ yielding  Algorithm~\ref{alg:det}. 
\begin{remark}[Estimator comparison] 
%The $\delta$-smoothing-based randomized difference constructs an estimator of $\nabla \widehat{f}^i_{t,\delta}(x)$ rather than $\nabla f_t^i(x)$ itself, while requiring only two function value queries without probing all $d$ dimensions. In contrast, the deterministic difference provides a direct estimate of $\nabla f_t^i(x)$ but requires $2d$ function value queries. We provide a more detailed comparison of their performance and practical trade-offs in Section~III.C after establishing the convergence guarantees of both algorithms.
The randomized difference estimates $\nabla \widehat{f}^i_{t,\delta}(x)$ rather than $\nabla f_t^i(x)$, but requires only two function-value queries. In contrast, the deterministic difference directly estimates $\nabla f_t^i(x)$ at the cost of $2d$ queries. Their theoretical performance and practical trade-offs are compared in Section~3.3.
\end{remark}
\begin{algorithm}[t]
\caption{Randomized Difference Based Algorithm}
\label{alg:rand}
\begin{algorithmic}[1]
\Require Same as in Framework~\ref{fr:common}.
\For{$t=0,\ldots,T$}
\State Each agent $i$ executes Steps 2--5 in Framework~\ref{fr:common} to obtain $S_t^i,z_{t+1}[i],y_t^i$.
\State \textbf{Randomized difference:} %Generate $u_t^i\sim \mathrm{Unif}(\mathbb S)$ independently of $\{u_s^j:~0\le s\le t-1,~j\in V\}$
Draw $u_t^i\sim\mathrm{Unif}(\mathbb S)$ independently of all information available before this draw, 
%$\{u_s^j\}_{0\le s\le t-1,~j\in V}$
and compute
\begin{align}
g_t^i\!=\!\frac{d}{2\delta_t}\left[f_t^i(y_t^i\!+\!\delta_t u_t^i)\!-\!f_t^i(y_t^i\!-\!\delta_t u_t^i)\right]u_t^i.\label{sa3}
\end{align}
\State \textbf{Update rule:} See (\ref{sa4}) in Framework~\ref{fr:common}. 
\EndFor
\end{algorithmic}
\end{algorithm}

\begin{algorithm}[t]
\caption{Deterministic Difference Based Algorithm}
\label{alg:det}
\begin{algorithmic}[1]
\Require Same as in Framework~\ref{fr:common}.
\For{$t=0,\ldots,T$}
\State Each agent $i$ executes Steps 2--5 in Framework~\ref{fr:common} to obtain $S_t^i,z_{t+1}[i],y_t^i$.
\State \textbf{Deterministic difference:} compute
\begin{align}
\nonumber h_t^i[k]&=\frac{f_t^i(y_t^i\!+ \!\delta_t\textbf{e}^k)\!\!-\!\!f_t^i(y_t^i\!-\!\delta_t \textbf{e}^k)}{2\delta_t},~k\!=\!1,\!\cdots\!,d, \\
h_t^i &= \left[h_t^i[1],\cdots,h_t^i[d]\right]^{\top}. \label{a3}
\end{align}
\State \textbf{Update rule:} See (\ref{sa4}) in Framework~\ref{fr:common}. 
\EndFor
\end{algorithmic}
\end{algorithm}

%remark：问题描述讲到拜占庭攻击时用图区别节点攻击与边攻击

\section{Consensus of Estimates and sublinearity of Regret}

%In this section, we establish two key theoretical properties of the proposed algorithms: (i) consensus of the agents’ estimates and (ii) bounds on the dynamic regret. A key step is to control the disagreement induced by the time-varying mixing matrices $\{\tilde{A}_t\}$ defined in (\ref{A}). We first establish a uniformly bound for ${z_t}$ together with an exponential mixing property for ${\tilde A_t}$, then we establish consensus and regret bounds for Algorithm~1, and finally extend the analysis to Algorithm~2. 

In this section, we establish consensus and dynamic regret bounds for the proposed algorithms. The analysis relies on controlling the disagreement induced by $\{\tilde A_t\}$. To this end, we first establish uniform bounds for $\{z_t\}$ and an exponential mixing property of $\{\tilde A_t\}$, then derive consensus and regret bounds for Algorithm~1, and finally extend the analysis to Algorithm~2. Define $$\tilde{\Phi}_{t,s}\triangleq \tilde{A}_{t}\tilde{A}_{t-1}\cdots \tilde{A}_{s+1}\tilde{A}_{s}\text{, }\forall t\geqslant s\text{, }\quad \tilde{\Phi}_{t,t+1}\triangleq I. $$
We summarize the required scaling and mixing properties in the following lemma.
\begin{lemma}\label{Lem:pi}
{Consider the two algorithms proposed in Section 2.3.} Assume Assumption~3 holds. Then,
\begin{itemize}
 \item[a)]There exist constants $\xi_1,\xi_2>0$, such that, for all $t\geqslant 0$ and all $i\in V$,
\begin{align}
\xi_1^{-1}\le z_t[i]\le \xi_2.\label{theta}
\end{align}
\item[b)] There exists a sequence of stochastic vectors $\{\pi_t\in \mathbb{R}^n,$ ${t} \geqslant 0\}$ and positive constants $C,~\beta$ and $\lambda\in(0,1)$ such that, for all $t\ge s\ge 0$ and all $i,j\in V$,
\begin{align}
   \mid\tilde{\Phi}_{t,s}[i,j]-\pi_s[j]\mid \le C\lambda^{t-s},~t\geqslant s\geqslant 0,\label{Phi'1}
\end{align}
with $\pi_t[i] \geqslant \beta$, and 
\begin{align}
\pi^{\top}_t=\pi^{\top}_{t+1} \tilde{A}_t,\quad \forall\, t\geqslant 0.\label{Phi'2}
\end{align}
\end{itemize}
\end{lemma}
\proof See \ref{Proof:Lem pi}. 

Before proceeding further, we need the following assumption for the perturbation radius. 
\begin{assumption}\label{asmp:delta}
$\{\delta_t\}_{t\geqslant 0}$ is a positive sequence nonincreasing to zero and  $\delta_0< r$ with $r$ defined in Assumption~\ref{asmp:X-interior}.
\end{assumption}
\subsection{Theoretical analysis for Algorithm 1}
Define $\overline{x}_{t}=\sum_{i=1}^n \pi_t[i] x_t^i$, where $\{x_{t}^i\}_{t\geqslant 0},i\in V$ are the estimates generated by Algorithm 1. Since $X$ is convex and $\pi_t$ is stochastic, we have $\bar x_t\in X$ for all $t$. 

\begin{theorem}[Consensus of Estimates]\label{sThm1}
Consider Algorithm \ref{alg:rand} defined in Section 2.3. Assume that Assumptions 1--4 hold.
\begin{itemize}
\item[i)] For any $i\in V$, it holds that
\begin{equation*}
\left\|x_{t}^i-\overline{x}_{t}\right\|\leq{\beta}_{t}, 
\end{equation*}
where
\begin{align}
\nonumber\beta_t\triangleq &
nCM\lambda^{t-2}+n\xi_1 C L_g\sum_{s=2}^{t-1}{\alpha_{s-1}}\lambda^{t-s-1}
\\&~+2\xi_1L_g\alpha_{t-1},\label{s86}
\end{align}
for $t\geqslant 3$ and $\beta_1=\beta_2\triangleq 2M$. Here we define $L_g\triangleq dL$, and $C>0,\lambda\in(0,1)$ are the constants in \eqref{Phi'1}.

\item[ii)] If $\alpha_{t}\rightarrow 0 $ as $t\rightarrow\infty$, then ${\beta}_{t}\rightarrow 0 $, i.e., Algorithm 1 achieves the consensus.
\end{itemize}
\end{theorem}
\proof See \ref{Proof:Thm 1.1}. 

%\begin{remark}
%Theorem~\ref{sThm1} characterizes the disagreement between each local estimate $x_t^i$ and the weighted network average $\bar{x}_t$. Specifically, part~i) provides an explicit upper bound $\beta_t$ for the estimation error $\|x_t^i-\bar{x}_t\|$, which quantifies how far the agents are from consensus at time $t$. Moreover, part~ii) shows that when the stepsize sequence $\{\alpha_t\}$ decreases to zero, the bound $\beta_t$ also converges to zero. This implies that the disagreement among the agents vanishes asymptotically, and hence Algorithm~1 achieves consensus.

%Theorem~\ref{sThm1} provides an explicit disagreement bound and implies asymptotic consensus for Algorithm~1.
%\end{remark}

Before bounding the dynamic regret, we establish a one-step descent inequality. Define $\theta_t\triangleq\|x_t^*-x_{t+1}^*\|$ and $\tilde{x}_t^*\triangleq (1-\delta_t/r)x_t^*$ for $t\ge 0$. 

%for $\|\bar x_t-x_t^*\|^2$. 
%to establish the gap between 
\begin{lemma}\label{sLem4.2}
Assume Assumptions 1--4 hold and $\{\alpha_t\}_{t\geqslant 0}$ is a positive sequence. Let $\{x_t^i\}_{t\geqslant 1,~i\in V}$ be the estimates generated by Algorithm 1. If for any $i\in V\text{ and }1\le t \le T$, $f_t^i(\cdot)$ is $\sigma$-strongly convex on $X$, then 
\begin{align}
\nonumber&\frac{1}{2}\mathbb{E}\left\|\bar{x}_{t\!+\!1}\!-\!\tilde{x}_{t\!+\!1}^*\right\|^2\!-\!\frac{1}{2}\mathbb{E}\left\|\bar{x}_t\!-\!\tilde{x}_t^*\right\|^2
\leq 2M\theta_t\!+\!\frac{1}{2}\xi_1^2\alpha_t^2L_g^2
\\\nonumber&+2\xi_1L_g\alpha_t\beta_{t\!+\!1}+{n\left(2L\!+\!LM/r\right)\alpha_t\delta_t}\!{+\!{2M^2}/{r}(\delta_t\!-\!\delta_{t+1})}\\
\nonumber&+\sum_{i=1}^n\alpha_t\mathbb{E}\left[f_t^i\left(x_t^{*}\right)-f_t^i\left(\bar{x}_t\right)-\frac{\sigma}{2}\left\|\bar{x}_t-\tilde{x}_t^*\right\|^2\right]
\\&+\!n\left(2M\sigma\!\!+\!L\!+\!\!L_g\right)\alpha_t\beta_t. \label{s8}
\end{align}
\end{lemma}
\proof See \ref{Proof:Lem 1.1}. 

Next we analyze the upper bound of the dynamic regret  (\ref{5}) for Algorithm 1. We first define a set of constants to be used later on:
\begin{equation*}%\label{s10}
\begin{split}
C_1&=\left(\!{2\xi_1 L_g}\!+n(2M\sigma\!+\!2 L+L_g)\right)\!\left(\!\beta_1\!\!+\!\!\beta_2\!\!+\!\!\frac{n C M}{1\!-\!\lambda}\!\right),\\
C_2&=\!\frac{\xi_1^2 L_g^2}{2 }\!+\!\left(\!{2\xi_1 L_g}\!+\!n(2M\sigma\!+\!2 L\!+\!L_g)\right)\!\left(\!\frac{n\xi_1 C L_g}{1-\lambda}\!+\!2\xi_1 L_g\!\right),\\
C_3&=\!\frac{2M^2}{ \alpha_1}+C_1,~{C_4=n(2L\!+\!LM/r)}.
\end{split}
\end{equation*}
where $L_g=dL$, $\beta_1=\beta_2=2M$, and  $M,~L,~\xi_1,~C,~\lambda,~r$ are constants specified in Assumptions  1--4, (\ref{theta}), (\ref{Phi'1}). 
\begin{theorem}\label{sThm2}
Assume that Assumptions 1--4 hold and $\{\alpha_t\}_{t\geqslant 1}$ is a positive decreasing sequence. Let $\{x_t^i\}_{t\geqslant 1,~i\in V}$ be the estimates generated by Algorithm 1. 
\begin{itemize}
\item[i)] If for $i\in V$ and $0\le t \le T$, $f_t^i(\cdot)$ is convex on $X$, then
\begin{align}
\nonumber\mathbb{E}\left[\mathrm{regret}_T^{D,i}\right] &\le C_1+\frac{2M^2}{\alpha_T}+C_2\sum_{t=1}^T \alpha_t 
+{2 M}\sum_{t=1}^T \frac{\theta_t}{\alpha_t}
\\&~+{{C_4}\sum_{t=1}^T\delta_t+\frac{2M^2}{r}\sum_{t=1}^T\frac{\delta_t\!-\!\delta_{t+1}}{\alpha_t}}.\label{s18a}
\end{align}
Moreover, if choosing $\alpha_t=1/\sqrt{t+1}$,  $\delta_t=C_{\delta}/\sqrt{t+1}$ with $0<C_{\delta }< r$,  and $\theta_t/\alpha_t=o(1)$, it holds that
\begin{align}
%\nonumber&
\mathbb{E}\left[\mathrm{regret}_T^{D,i}\right]%\le C_1 +(2M^2\!+\!C_2\!+\!{C_4\!+\!\frac{2M^2}{r}})O(\sqrt{T})\\&~+2M\sum_{t=1}^T \frac{\theta_t}{\alpha_t}
=O(\sqrt{T})+o(T). \label{s18aa}
\end{align}
\item[ii)] If for any $i\in V\text{ and }0\le t \le T$, $f_t^i(\cdot)$ is $\sigma$-strongly convex on $X$, then for $i\in V$
\begin{align}
\nonumber&\mathbb{E}\left[\mathrm{regret}_T^{D,i}\right]  \le C_3 +C_2\sum_{t=1}^T \alpha_t +{2 M}\sum_{t=1}^T \frac{\theta_t}{\alpha_t}\\
\nonumber&+{{{C_4}}\sum_{t=1}^T\delta_t +\frac{2M^2}{r}\sum_{t=1}^T\frac{\delta_t\!-\!\delta_{t+1}}{\alpha_t}}
\\&+\frac{1}{2}\sum_{t=2}^T {\mathbb{E}\left\|\bar{x}_t-\tilde{x}_t^*\right\|^2}\left(\frac{1}{\alpha_t}-\frac{1}{\alpha_{t-1}}-\sigma n\right)\label{s18}
\end{align}
Moreover, if choosing $\alpha_t=1/(\sigma n (t+C_\alpha))$,  $\delta_t=C_{\delta}/(t+1)$ with $C_\alpha>0,~0<C_{\delta }< r$, and $\theta_t/\alpha_t=o(1)$,
\begin{align}
%\nonumber&
\mathbb{E}\left[\mathrm{regret}_T^{D,i}\right]
%\le C_3 +(C_2+{C_4}+2M^2/r)O(\log(T))\\&~+2M\sum_{t=1}^T \frac{\theta_t}{\alpha_t}
=O(\log(T))%+O(1)
+o(T). \label{s18b}
\end{align}
\end{itemize}
\end{theorem}
\proof See \ref{Proof:Thm 1.2}. 

%\begin{remark}
%Theorem~\ref{sThm2} provides upper bounds on the dynamic regret of Algorithm~1 under different assumptions on the objective functions, and further specifies the corresponding regret orders under particular choices of the stepsize and perturbation radius. 
%In the convex case, by choosing $\alpha_t=1/\sqrt{t}$ and $\delta_t=1/\sqrt{t}$, the dynamic regret is bounded by $O(\sqrt{T})+O(1)+o(T)$. In the strongly convex case, by choosing $\alpha_t=1/(\sigma n t)$ and $\delta_t=1/t$, the dynamic regret is bounded by $O(\log T)+O(1)+o(T)$. 
%A comparison between these regret orders and those reported in related works will be provided later after establishing the dynamic regret bound for Algorithm~2.

%Theorem~\ref{sThm2} establishes dynamic regret bounds for Algorithm~1 under different assumptions on the objective functions, and further specifies the corresponding sublinear regret orders under particular choices of the stepsize and perturbation radius. A comparison with related results will be given after presenting the corresponding bounds for Algorithm~2.
%\end{remark}

\subsection{Theoretical analysis for Algorithm 2}
To control the bias of the deterministic finite-difference estimator, we impose the following condition. 
\begin{assumption}
$\{f_t^i(\cdot), t\ge 0 \},~i\in V$ are second order differentiable on $X$ and there exists  $L_{H}>0$ such that
\begin{equation*}
\begin{split}
\left|{\nabla^2f_t^i}[k,k](x)\right|\le L_{H}
\end{split}
\end{equation*}
for any $x\in X~,i\in V$, $t\ge 0$ and $k\in \{1,...,d\} $, where $\nabla^2 f^i_t(\cdot)$ represents the Hessian matrix of $f^i_t(\cdot)$.
\end{assumption}
%%%根据online书Lemma2.6，有$\left|{\nabla^2f_t^i}[k,k]\right|\le L^{'}(x\in X)$,但是此处需要在$X'$中成立。

%Define 
%\begin{align*}\label{bar_x}
%\overline{x}_{t}=\sum_{i=1}^n \pi_t[i] x_t^i, 
%\end{align*}
%where $\{x_{t}^i\}_{t\geqslant 0},i\in V$ are the estimates generated by Algorithm 2. 
By carrying out a similar proof as that of Theorem \ref{sThm1}, we can obtain the following consensus result.
\begin{theorem}[Consensus of Estimates]\label{Thm1}
Assume that Assumptions 1--5 hold.
\begin{itemize}
\item[i)] For the estimates generated by Algorithm 2, it holds
\begin{equation}\label{T1}
\left\|x_{t}^i-\overline{x}_{t}\right\|\leq{\beta}_{t},
\end{equation}
for any $i\in V$, where $\overline{x}_{t}=\sum_{i=1}^n \pi_t[i] x_t^i$, and 
\begin{align}
\nonumber\beta_t & \triangleq
nC\max_{j\in V}\left\|x^{j}_{1}\right\|\lambda^{t-2}+n\xi_1 C L_h\sum_{s=2}^{t-1}{\alpha_{s-1}}\lambda^{t-s-1}\\
&\quad +2\xi_1L_h\alpha_{t-1},\label{86}
\end{align}
for all $t\geqslant 3$ and $\beta_1=\beta_2\triangleq 2M$, where $L_h\triangleq L+\delta_0\sqrt{d}L_H$, $C>0$ and $\lambda\in(0,1)$ are the constants specified in formula (\ref{Phi'1}).

\item[ii)] If $\alpha_{t}\rightarrow 0 $ as $t\rightarrow\infty$, then ${\beta}_{t}\rightarrow 0 $, i.e., Algorithm 2 achieves the consensus.
\end{itemize}
\end{theorem}
We now establish the regret bounds of the deterministic-difference based algorithm. 
\begin{equation*}%\label{10}
\begin{split}
C^{'}_1&\!=\!\left({2\xi_1 L_h}\!+n\left(2M\sigma+3L\right)\right)\!\left(\beta_1\!+\!\beta_2\!+\!\frac{nCM}{1\!-\!\lambda}\right),\\
C^{'}_2&\!=\!\frac{\xi_1^2\! L_h^2}{2 }\!\!+\!\!\left({2\xi_1 L_h}\!+n\left(2M\sigma+3L\right)\right)\left(\frac{n\xi_1 C\! L_h}{1\!-\!\lambda}\!\!+\!\!2\xi_1 L_h\!\right),\\
C^{'}_3&\!=\frac{2M^2}{ \alpha_1}+C^{'}_1,~{C^{'}_4=nM\left(\frac{L_h}{r}+2\sqrt d\,L_H\right)}.
\end{split}
\end{equation*}
where $L_h\!=\!L\!+\!\delta_0\!\sqrt{d} L_H,M,L,\delta_0,L_H,\xi_1,C,\lambda$ are constants specified in Assumptions 1,~2,~4,~5 and (\ref{theta}) -- (\ref{Phi'1}). 
\begin{theorem}\label{Thm2}
Assume that Assumptions 1--5 hold and $\{\alpha_t\}_{t\geqslant 1}$ is a positive decreasing sequence. Let $\{x_t^i\}_{t\geqslant 1,~i\in V}$ be the estimates generated by Algorithm 2. 
\begin{itemize}
\item[i)] If for $i\in V$ and $0\le t \le T$, $f_t^i(\cdot)$ is convex on $X$, then
\begin{align}
\nonumber\mathrm{regret}_T^{D,i}  \le& C^{'}_1+\frac{2M^2}{\alpha_T} +C^{'}_2\sum_{t=1}^T \alpha_t +{2 M}\sum_{t=1}^T \frac{\theta_t}{\alpha_t}
\\&+{C^{'}_4} \sum_{t=1}^T \delta_t,\label{18a}
\end{align}
Moreover, if choosing $\alpha_t=1/\sqrt{t+1}$,  $\delta_t=C_{\delta}/\sqrt{t+1}$ with $0<C_{\delta }< r$, and $\theta_t/\alpha_t=o(1)$, it holds that
\begin{align}
%\nonumber&
\mathrm{regret}_T^{D,i}
%\le C_1^{'} +(2M^2+C_2^{'}+C_4^{'})O(\sqrt{T})\\&~+2M\sum_{t=1}^T \frac{\theta_t}{\alpha_t}
=O(\sqrt{T})%+O(1)
+o(T). \label{18aa}
\end{align}
\item[ii)] If for any $i\in V\text{ and }0\le t \le T$, $f_t^i(\cdot)$ is $\sigma$-strongly convex on $X$, then for $i\in V$
\begin{align}
\nonumber&\mathrm{regret}_T^{D,i}  \le C^{'}_3 +C^{'}_2\sum_{t=1}^T \alpha_t +{2 M}\sum_{t=1}^T \frac{\theta_t}{\alpha_t}+{C^{'}_4}\sum_{t=1}^T \delta_t
\\&~~~~+\frac{1}{2}\sum_{t=2}^T {\left\|\bar{x}_t-x_t^*\right\|^2}\left(\frac{1}{\alpha_t}-\frac{1}{\alpha_{t-1}}-\sigma n\right),\label{18}
\end{align}
Moreover, if choosing $\alpha_t=1/(\sigma n (t+C_\alpha))$,  $\delta_t=C_{\delta}/(t+1)$ with $C_\alpha>0,~0<C_{\delta }< r$, and $\theta_t/\alpha_t=o(1)$, 
%and terms $\sum_{t=1}^T{c_t})$, $\sum_{t=1}^T\frac{\theta_t}{\alpha_t})$ and $\sum_{t=1}^T{\tilde{\theta}_t}$ are bounded by $o(T)$, 
\begin{align}
%\nonumber&
\mathrm{regret}_T^{D,i}
%\le C_3^{'} +(C_2^{'}+C_4^{'})O(\log(T)) \\&~+2M\sum_{t=1}^T \frac{\theta_t}{\alpha_t}
=O(\log(T))%+O(1)
+o(T). \label{18b}
\end{align}
\end{itemize}
\end{theorem}

{\emph{Proof of Theorem \ref{Thm1}--\ref{Thm2}: }} See \ref{Proof:Thm2}. 

{
\begin{remark} [Online versus static optimization] 
Unlike static optimization, online optimization requires agents to act before observing the current loss, which may vary across rounds. Dynamic regret measures cumulative excess loss relative to each round's minimum, so its analysis accounts for minimizer movement, disagreement among agents, and zeroth-order estimation errors. For time-invariant objectives, a fixed minimizer gives zero path variation. Under suitable conditions, Algorithm~1 converges almost surely by the Robbins-Siegmund theorem, and Algorithm~2 converges deterministically. Under strong convexity and suitable parameter choices, the squared distance to the minimizer is $O(1/t)$ in expectation for Algorithm~1 and deterministically for Algorithm~2. %The technical report~\cite{liu2025zeroth-arxiv} gives the full statements and proofs.
\end{remark}}

{\begin{corollary}[Time-invariant objectives]\label{cor:static}
Suppose that the convex local objectives $f_t^i$ do not depend on $t$, i.e., $f_t^i=f^i$ for $t\ge 0$, and that their sum $f=\sum_{i=1}^n f^i$ has a unique minimizer $x^*$ on $X$. Assume that Assumptions~1--4 hold for Algorithm~\ref{alg:rand} and Assumptions~1--5 hold for Algorithm~\ref{alg:det}.
\begin{itemize}
\item[i)] If the positive stepsizes satisfy
\[
\sum_{t=1}^{\infty}\alpha_t=\infty,\quad
\sum_{t=1}^{\infty}\alpha_t^2<\infty,\quad
\sum_{t=1}^{\infty}\alpha_t\delta_t<\infty,
\]
then for all $i\in V$, 
\[
\begin{aligned}
&\lim_{t\to \infty}x_t^i= x^* &&\text{almost surely for Algorithm ~\ref{alg:rand}},\\
&\lim_{t\to \infty}x_t^i= x^* &&\text{deterministically for Algorithm~\ref{alg:det}}.
\end{aligned}
\]
\item[ii)] If each $f^i$ is $\sigma$-strongly convex with $\sigma>0$, choose $\alpha_t=1/[n\sigma(t+1)]$ and $\delta_t=C_\delta/(t+1)$ with $0<C_\delta<r$. Then, for all $i\in V$,
\[
\begin{aligned}
\mathbb E\|x_t^i-x^*\|^2&=O(1/t) &&\text{for Algorithm~\ref{alg:rand}},\\
\|x_t^i-x^*\|^2&=O(1/t) &&\text{for Algorithm~\ref{alg:det}}.
\end{aligned}
\]
\end{itemize}
\end{corollary}
\proof See \ref{Proof:static}.}

\begin{remark}[Path variation] %[Comparison for dynamic regret bounds between this work and existing results without attacks]
%{We next compare our dynamic regret bounds with existing results for online optimization without Byzantine edge attacks.} 
Dynamic regret bounds typically involve a sum of the variation term $\|x_{t+1}^*-x_t^*\|$, which quantifies the temporal regularity of the per-round minimizers. 
Intuitively, no online algorithm can track $x_t^*$ well when the minimizer varies significantly %from one round to the next
~\cite{2018Shahrampour,2022Xiong,2019Zhang}. 
For example, the dynamic regret bound in~\cite{2018Shahrampour} and \cite{2022Xiong} can be interpreted as having the order 
$O(\sqrt{T}) + O\!\left(\sum_{t=1}^T \|x_{t+1}^*-x_t^*\|/\alpha_t\right)$, in the gradient-based and zeroth-order convex setting, respectively. Moreover, 
\cite{2019Zhang} derives a gradient-based bound of order 
$O\!\Big(1+\sum_{t=1}^T \|x_{t+1}^*-x_t^*\|+\sum_{t=1}^T \|g_{t+1}-g_t\|_{\infty}\Big)$ with $g_t[i]\triangleq \nabla f_t^i(x_t^i)$ for $i=1,\ldots,n$. Theorems~\ref{sThm2} and~\ref{Thm2} show that, under Byzantine edge attacks, our results 
%two distributed zeroth-order algorithms attain dynamic regret bounds of order $O(\sqrt{T}) + O\!\left(\sum_{t=1}^T \|x_{t+1}^*-x_t^*\|/\alpha_t\right)$ for convex objectives, and $O(\log T) + O\!\left(\sum_{t=1}^T \|x_{t+1}^*-x_t^*\|/\alpha_t\right)$ for $\sigma$-strongly convex objectives. These rates 
are comparable to the dynamic regret guarantees without attacks.
\end{remark}

\subsection{Comparison between the two algorithms}
%Algorithm~1 adopts a randomized difference scheme, which requires only two function value queries per iteration and does not impose smoothness assumptions on the objective functions. These features make Algorithm~1 applicable to nonsmooth problems or scenarios where function-value sampling is expensive. However, its convergence guarantees are established only in expectation, and the convergence rate deteriorates as the problem dimension increases. 

Algorithm~1 uses a randomized difference scheme, requiring only two function-value queries per iteration, without the additional second-order smoothness conditions in Assumption~5. It is therefore suitable when sampling is limited or these smoothness conditions cannot be assumed. However, its sublinear regret guarantees hold only in expectation, and its performance deteriorates as the dimension increases.

%Algorithm~1 uses a randomized difference scheme requiring only two function-value queries per agent per iteration. Unlike Algorithm~2, it does not require the additional second-order smoothness condition in Assumption~5. Its consensus guarantee holds pathwise, whereas its dynamic regret bound is established in expectation. Algorithm~2 uses deterministic coordinate-wise differences, requiring $2d$ function-value queries per agent per iteration, and provides pathwise regret bounds under the additional smoothness condition. 

%Specifically, by $L_g=dL$, $L_h=L+\delta_0\sqrt{d}L_H$, and the definitions of the constants $C_i$ and $C_i'$, $i=1,2,3$, for dimension $d$,it can be readily verified that
Specifically, with all other parameters fixed independently of $d$,
the identities $L_g=dL$, $L_h=L+\delta_0\sqrt{d}L_H$,
and the definitions of $C_i$ and $C_i'$, $i=1,\ldots,4$, yield
\begin{align*}
C_2&=O(d^2),~C_1=C_3=O(d),~C_4=O(1),\\
C_2^{'}&=O(d),~C_1^{'}=C_3^{'}=O(\sqrt{d}),~C_4^{'}=O(\sqrt{d}).
\end{align*}
%\begin{align*}
%\nonumber0&\le C_i-C_i^{'}=O\left(3L_g-2L_h-L\right)
%\\&=O\left(3Ld-2\delta_1L_H\sqrt{d}\right), 
%\end{align*}
%holds for $i=1,3$ and 
%\begin{align*}
%\nonumber&0\le C_2-C_2^{'}=O(L_g^2-LL_h)+O\left(L_g-L_h\right)
%\\=&O\left(L^2d^2\!-\!\delta_1L_HL\sqrt{d}\right)+O\left(Ld\!-\!\delta_1L_H\sqrt{d}\right), 
%\end{align*}
%holds for $d\in \{d:Ld-\delta_1L_H\sqrt{d}-L\geqslant 0\}$. }
%These relations imply that, as the dimension $d$ increases, the convergence constants and regret bounds of Algorithm~1 become significantly larger than those of Algorithm~2. This trend is consistent with Theorems~1--4 and is further corroborated by the numerical results in Section~IV.
%In summary, Algorithm~1 is preferable in applications with limited sampling capability or nonsmooth objectives, whereas Algorithm~2 is more suitable for scenarios demanding deterministic performance guarantees or fast convergence. 
In addition, the regret bounds of Algorithm~1 contain a term $\frac{2M^2}{r}\sum_{t=1}^T\frac{\delta_t\!-\!\delta_{t+1}}{\alpha_t}$. 
%As $d$ increases, the convergence constants and regret bounds of Algorithm~1 become significantly larger than those of Algorithm~2, which indicates milder explicit dimension dependence for Algorithm~2. This is consistent with Theorems~1--4 and is also confirmed by the numerical results in Section~4.
With the remaining parameters held dimension-independent, the leading dimension-dependent terms grow more rapidly with $d$ for Algorithm~1 than for Algorithm~2, suggesting Algorithm~2 may exhibit more favorable regret performance in higher-dimensional problems. %This is also confirmed by the numerical results in Section~4. 
The numerical results in Section 4 provide empirical evidence of this qualitative trend under the considered simulation settings.

In summary, Algorithm~1 is preferable when sampling is limited or second-order smoothness cannot be assumed, whereas Algorithm~2 is more suitable when deterministic guarantees or fast convergence are required. 
{
\begin{remark}[Zeroth-order and first-order]
%0阶很重要，不只是一阶变为零阶
%理论上有什么差别
%一阶后也work

Zeroth-order methods use function values to estimate gradients in black-box problems~\cite{2005Flaxman,2012Shalevshwartz}, whereas first-order methods use gradient feedback directly. In our framework, the zeroth-order analysis accounts for smoothing and random estimation errors in Algorithm~1, finite-difference bias in Algorithm~2, and domain-contraction errors in both. With exact first-order feedback, both algorithms use the same resilient mixing rule, set
$G_t^i=\nabla f_t^i(y_t^i)$, and use $X$ for projection and for the feasibility check in $S_t^i$. These choices eliminate the approximation
and contraction errors. Extending the analysis to stochastic first-order feedback would require conditional unbiasedness, bounded conditional variance, and bounds on the effects of gradient noise on disagreement and regret.
\end{remark}}

\section{Simulation}
%We consider a multi-agent system with $n$ agents communicating over a stochastically generated time-varying directed network. Each agent has a trusted self-loop. During the first $B_0=5$ rounds, each possible nonself directed edge is independently generated with probability $q_n=\min\{100,n\}/n$ and classified as trusted, adversarial, or normal with probabilities $0.2$, $0.1$, and $0.7$, respectively. If the union of the trusted edges is not strongly connected, additional trusted edges are added or existing nontrusted edges are reclassified in the last round to ensure strong connectivity, based on the strongly connected components of the union graph. The resulting sequence of graphs is then repeated periodically, ensuring that the union of trusted edges over every $B_0$ consecutive rounds is strongly connected.
{
In this section, we numerically validate the theoretical results and compare convergence speed and dynamic regret across the proposed and baseline algorithms.

We consider a system with $n$ agents over a stochastically generated directed network with time-varying edge types. Each agent has a trusted self-loop. During the first $B_0=5$ rounds, all possible nonself directed edges are present, and each edge is independently 
%generated with probability $q_n=\min\{100,n\}/n$ and 
classified as trusted, adversarial, or normal with probabilities $0.2$, $0.03$, and $0.77$, 
%$0.1$, and $0.7$, 
respectively. If necessary, trusted edges are added or nontrusted edges are reclassified in the last round to ensure that their union is strongly connected. The resulting graph sequence is repeated periodically.

We consider a tracking model widely used in aircraft tracking, missile interception, etc~\cite{2020Yi}. The loss function of agent $i$ at time $t$ is defined by
\begin{equation*}
\begin{aligned}
f_{ t}^i\left(x\right) &=\zeta_{i, 1}\left\langle\omega^i, x\right\rangle+\zeta_{i, 2}\left\|x-\xi_{t}\right\|^{2}.
\end{aligned}
\end{equation*}
Assume $\zeta_{i, 1}=\zeta_{i, 2}=1$ for $i\in V$. Choose $X=[-5,5]^d\subset \mathbb{R}^d$, $d=4$, $\xi_t[k]=\frac{2d-k}{d(t+1)},~k=1, \ldots, d$, $\omega^i$ to be i.i.d. with uniform distribution over $[0,1]^d$. Each coordinate of $x_0^i$ is independently sampled from the uniform distribution over $[-2,2]$. Choose $\alpha_t=\frac{1}{2n(t+1)}$, $\delta_t=\frac{1}{n(t+1)}$ for Algorithms 1--2. The attack information is modeled as a $d$-dimensional random vector whose entries are i.i.d. and follow the Gaussian distribution $\mathcal{N}(1,5^2)$. 

%%%%%%%%%为了测试算法1,2的动态遗憾的次线性以及对攻击的鲁棒性
%To evaluate the dynamic regret of Algorithms~1 and~2, for $n=50$, we consider the following four baselines: the adapted RDO-T algorithm, a trusted-skeleton variant of Algorithm~1, and two attack-free zeroth-order DOCO algorithms~\cite{2022Xiong,2023Liu}. The adapted RDO-T modifies the node-resilient RDO-T in~\cite{2019Zhao} by replacing its scalar filtering rule and first-order update with the multidimensional set $S_t^i$ and a randomized difference update at $x_t^i$, respectively. %, while retaining neighbor averaging without convex-hull reconstruction or push-sum normalization. 
%The trusted-skeleton algorithm retains the estimator and push-sum updates of Algorithm~1 but uses only trusted edges. The two attack-free baselines use the parameters specified in \cite{2022Xiong,2023Liu}. For the zeroth-order push-sum DOCO algorithm~\cite{2022Xiong}, only state messages are subject to attacks, while normalization scalars remain uncorrupted. Fig.~1 compares the time-averaged dynamic regret of the six algorithms, measured by $\max_{i\in V}\mathrm{regret}_t^{D,i}/t$. 
To evaluate the dynamic-regret performance of Algorithms~1 and~2, we consider four baselines with \(n=50\). First, to assess the applicability of node-resilient methods to Byzantine edge attacks, we include an adapted RDO-T algorithm based on the node-resilient RDO-T in~\cite{2019Zhao}, replacing its scalar filtering rule and first-order update with the multidimensional resilient neighbor set \(S_t^i\) and a randomized zeroth-order difference estimator evaluated at \(x_t^i\), respectively. Second, the trusted-skeleton algorithm retains the estimator and push-sum updates of Algorithm~1 but restricts communication to trusted edges. Finally, we include two attack-free zeroth-order DOCO algorithms~\cite{2022Xiong,2023Liu}, using the parameter settings specified in ~\cite{2022Xiong,2023Liu}. For the zeroth-order push-sum DOCO algorithm~\cite{2022Xiong}, only state messages are subject to attacks, while normalization scalars remain uncorrupted. With a logarithmic vertical axis, Fig.~1 compares the time-averaged dynamic regret $\max_{i\in V}\mathrm{regret}_t^{D,i}/t$ of all six algorithms.
%, averaged over $10$ paired runs.
\begin{figure}[t]
\begin{center}
\includegraphics[width=0.95\linewidth]{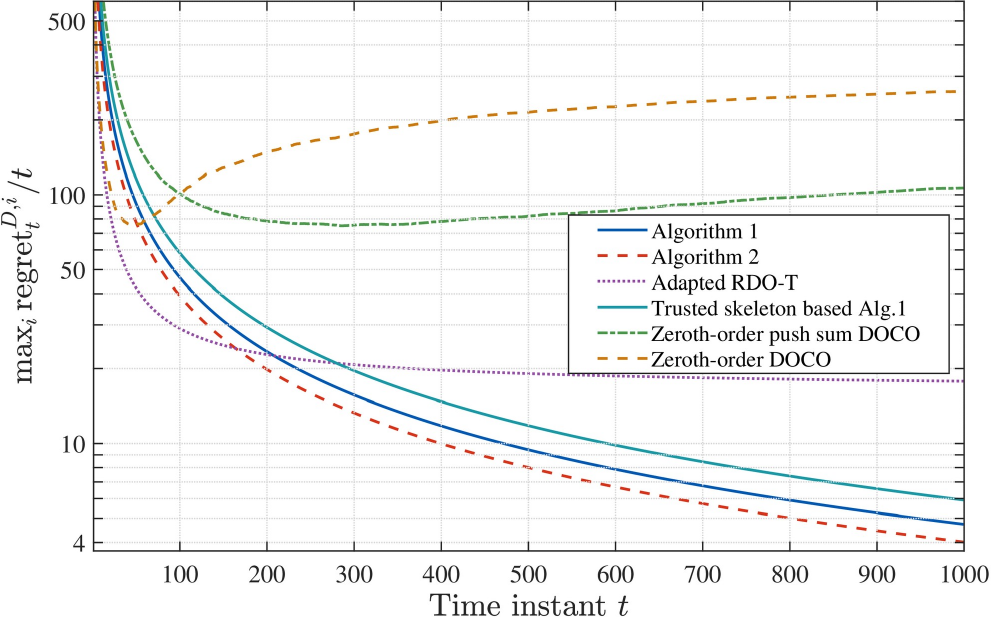}    
\caption{{Values of $\max_{i\in V}\mathrm{regret}^{D,i}_t/t$ of Algorithm 1, 2 in this work, the adapted RDO-T (proposed for node attacks in \cite{2019Zhao}, zeroth and high-dimension version), trusted skeleton (using only trusted edges), zeroth-order push sum DOCO (designed for non-attack scenarios in \cite{2022Xiong}), and zeroth-order DOCO algorithm (designed for non-attack scenarios in \cite{2023Liu}).} }  %
\label{fig3}                                 % 
\end{center}                                 % 
\end{figure}

To further compare Algorithms 1--2, Fig.~2 plots the consensus error 
$\max_{i\in V}\left\|x_t^i\!-\!{\sum_{j=1}^n x_t^j}/{n}\right\|$ with $n=50$, 
%To examine the effect of the dimension, we set  %to reduce the computational cost 
while Fig.~3 plots the number of rounds required to reach $\max_{i \in V} \mathrm{regret}^{D,i}_T / T \le 1$ for different $d$ with $n=10$.

%\begin{figure}[t]\begin{center}
%\includegraphics[width=8.6cm,height=3.0cm]{consensus2.eps}   \caption{{$\max_{i\in V}\left\|x_t^i-{\sum_{i=1}^n x_t^i}/{n}\right\|$ of Algorithms 1,~2, validating the consensus properties proposed in Theorems \ref{sThm1}, \ref{Thm1}. }}  \label{fig2}                       
%\end{center}\end{figure}

\begin{figure}[t]
\centering
\begin{minipage}[t]{0.48\linewidth}
    \centering
    \includegraphics[width=\linewidth]{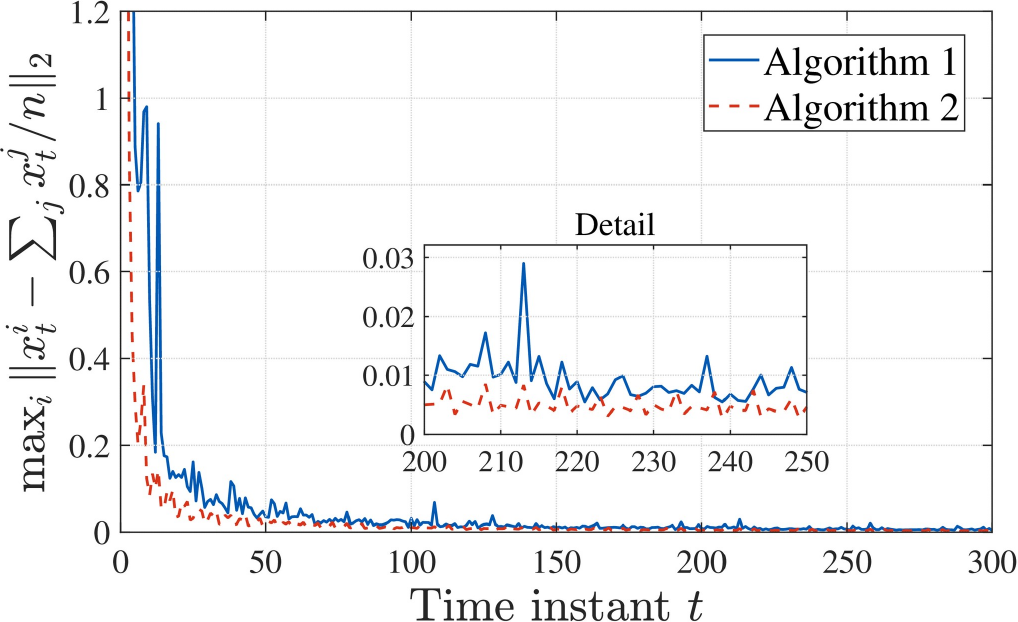}\caption{{
Consensus error
$\max_{i}\left\|x_t^i-
\frac{1}{n}\sum_{j=1}^n x_t^j\right\|$
of Algorithms~1 and 2, validating the consensus properties in Theorems~\ref{sThm1} and~\ref{Thm1}. }}
\end{minipage}
\hfill
\begin{minipage}[t]{0.48\linewidth}
    \centering
    \includegraphics[width=\linewidth]{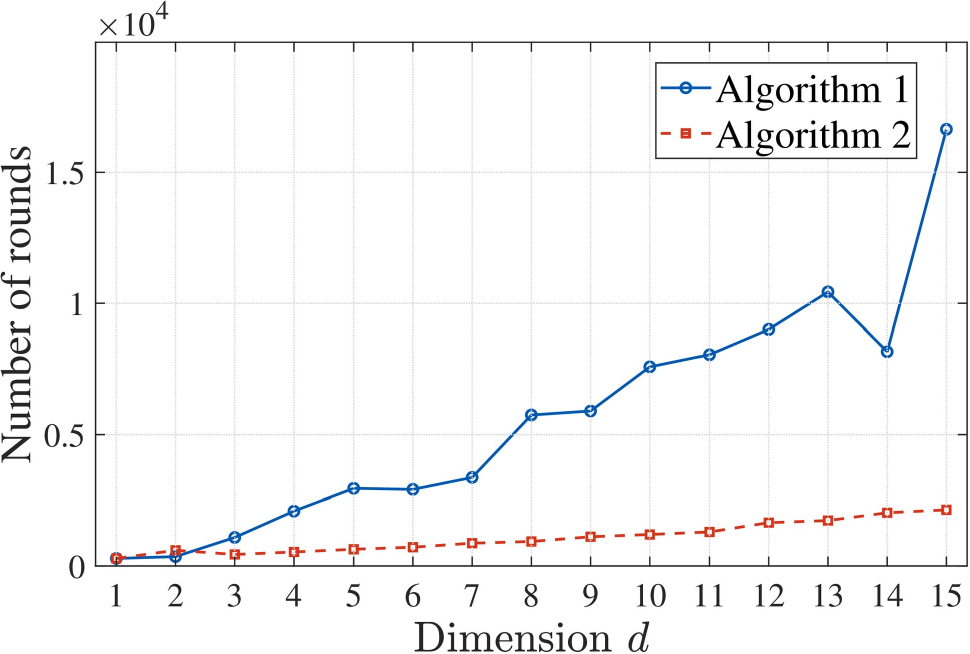}\caption{{
Number of rounds required by Algorithms~1 and~2
to achieve $\max_{i\in V}\mathrm{regret}^{D,i}_T/T
\le 1$ as $d$ varies, over a newly generated network
with $10$ agents.
}}
\end{minipage}
\label{fig:simulation}
\end{figure}

The preceding random-network experiments with i.i.d. local objective coefficients show modest gains over the trusted-skeleton (TS) baselines. We further consider a complete
directed network with $n=2m$ agents partitioned into $A=\{1,\ldots,m\}$ and $B=\{m+1,\ldots,2m\}$. All intra-cluster links and self-loops are trusted. The bridge $m\leftrightarrow m+1$ is trusted when $t\equiv1\pmod{20}$ and normal otherwise. The links $i\leftrightarrow i+m$, $i\in A$, are persistently attacked; all remaining links are normal. Let $v=\mathbf{1}_d/\sqrt d$, with $s_i=1$ on $A$ and $s_i=-1$ on $B$. Set $\omega^i=-2c_i$ and $x_0^i=2s_iv+\eta_i$, where $c_i=s_iv+\epsilon_i$, except that $c_m=2v+\epsilon_m$ and $c_{m+1}=-2v+\epsilon_{m+1}$.
The entries of $\epsilon_i$ and $\eta_i$ are independently sampled from $\mathrm{Unif}[-0.05,0.05]$ once per run. We use $n=10$ and $d=2$; other parameters remain unchanged. 
TS.~1 and TS.~2 differ from Algorithms~1 and~2 only by removing nontrusted-message contributions. 
Table~2 reports mean stable-attainment rounds over $10$ paired runs for the regret and consensus thresholds, where a criterion must remain satisfied from attainment through the horizon. Setting~I uses $B_0=20$ and $T=20000$. In Setting~II, successive trusted-bridge intervals instead grow as $20,30,40,\ldots$, with $T=500000$. This unbounded interval sequence violates Assumption~3 and tests performance beyond the established theoretical conditions.
%\begin{table}[t]\centering
%\caption{{Median number of rounds to reach and maintain the prescribed thresholds over $10$ paired runs.}}\label{tab:iterations}
%\begingroup\scriptsize\setlength{\tabcolsep}{2pt}\renewcommand{\arraystretch}{1.3}
%\begin{tabular}{|c|c|c|}
%\hline
%Criterion
%& Alg. 1/TS. 1
%& Alg. 2/TS. 2
%\\ \hline
%$\max_{i } \mathrm{regret}^{D,i}_t/t \le 5$
%& 585.5/1065
%& 390/1058.5
%\\ \hline
%$\max_{i } \mathrm{regret}^{D,i}_t/t \le 1$
%& 3427.5/6229
%& 2266/6232
%\\ \hline
%$\max_{i}\left\|x_t^i-{\sum_{i=1}^n x_t^i}/{n}\right\|\le 0.1$
%& 2223/3583
%& 1323/3573
%\\ \hline
%\end{tabular}\endgroup\end{table}

%\begin{table}[t]
%\centering \color{blue}
%\caption{{Median number of rounds to reach and maintain the thresholds for $r_t\!=\!\max_{i\in V} \mathrm{regret}^{D,i}_t/t$ and $e_t\!=\!\max_{i\in V}\left\|x_t^i\!-\!{\sum_{i=1}^n x_t^i}/{n}\right\|$ over $T\!=\!20000$ and $10$ paired runs.}}
%\label{tab:iterations}
%\begingroup
%\scriptsize
%\setlength{\tabcolsep}{2pt}
%\renewcommand{\arraystretch}{1.3}
%\begin{tabular}{|c|c|c|c|c|}
%\hline
%Criterion &$r_t \le 5$ &$r_t \le 1$ &$e_t\le 0.5$ & $e_t\le 0.1$
%\\ \hline
% Alg. 1/TS. 1 
% & 683/1130.5 %& 585.5/1065
% & 3905.5/6573.5  %& 3427.5/6229
% & 373/632.5 %& 343/623 
%& 2313/3602.5 %& 2223/3583
%\\ \hline
% Alg. 2/TS. 2 
% & 473/1106 %& 390/1058.5 
% & 2673.5/6479 %& 2266/6232  
% & 243/633 %& 223/623  
% & 1343/3593 %& 1323/3573 
% \\ \hline
%\end{tabular}
%\endgroup
%\end{table}

\begin{table*}[t]
\centering
%\color{blue}
\caption{{Mean number of rounds to reach and maintain the prescribed thresholds under two experimental settings.}}
\label{tab:iterations}
\begingroup
\fontsize{5pt}{5pt}\selectfont
\setlength{\tabcolsep}{1.2pt}
\renewcommand{\arraystretch}{1.3}
%\begin{tabular}{|c|ccccc|ccccc|ccccc|ccccc|}
\begin{tabular}{|c|c|c|c|c|c|c|c|c|c|c|c|c|c|c|c|c|c|c|c|c|}
\hline
& \multicolumn{10}{c|}{\textbf{Setting I: $T=20000$, $B_0=20$}}
& \multicolumn{10}{c|}{\textbf{Setting II: $T=500000$, growing bridge intervals}}
\\ \hline
Criterion
& \multicolumn{5}{c|}{$\max_{i} \mathrm{regret}^{D,i}_t/t \le r_0$}
& \multicolumn{5}{c|}{$\max_{i}\left\|x_t^i\!-\!{\sum_{j=1}^n\! x_t^j}/{n}\right\|\!\le \!r_0$}
& \multicolumn{5}{c|}{$\max_{i} \mathrm{regret}^{D,i}_t/t \le r_0$}
& \multicolumn{5}{c|}{$\max_{i}\left\|x_t^i\!-\!{\sum_{j=1}^n\! x_t^j}/{n}\right\|\!\le \!r_0$}
\\ \hline

$r_0$
& 5 & 4 & 3 & 2 & 1
& 0.5 & 0.4 & 0.3 & 0.2 & 0.1
& 5 & 4 & 3 & 2 & 1
& 0.5 & 0.4 & 0.3 & 0.2 & 0.1
\\ \hline

Alg.~1
& 675.6 & 877.5 & 1213.0 & 1882.8 & 3877.2
& 369.0 & 499.0 & 703.0 & 1103.0 & 2297.0
& 3560.9 & 5828.9 & 10230.4 & 20466.1 & 58687.4
& 5429.0 & 10097.0 & 19277.0 & 55451.0 & 253579.0
\\ \hline

TS.~1
& 1125.3 & 1473.4 & 2046.3 & 3182.6 & 6574.3
& 632.9 & 845.0 & 1180.8 & 1791.0 & 3600.8
& 10609.1 & 17143.6 & 29797.3 & 59648.1 & 170497.3
& 13669.0 & 25554.0 & 52731.0 & 133985.0 & NR
\\ \hline

Alg.~2
& 474.5 & 613.9 & 844.3 & 1303.5 & 2677.7
& 243.0 & 320.9 & 441.0 & 663.0 & 1345.0
& 1707.6 & 2605.8 & 4322.3 & 8321.5 & 22991.0
& 1993.0 & 3691.9 & 7403.0 & 18842.0 & 83070.0
\\ \hline

TS.~2
& 1106.0 & 1448.6 & 2013.3 & 3133.5 & 6479.0
& 633.0 & 842.9 & 1173.0 & 1791.0 & 3593.0
& 10543.6 & 17054.8 & 29668.8 & 59449.4 & 170233.0
& 13669.0 & 25698.0 & 52833.0 & 133985.0 & NR
\\ \hline
\end{tabular}

\vspace{1mm}
\parbox{\textwidth}{\tiny NR: no paired run reached and maintained the threshold within the horizon; numerical entries: all paired runs did.}
\endgroup
\end{table*}

Figs.~1 and~2 illustrate that Algorithms~1 and~2 achieve consensus and exhibit sublinear dynamic regret, consistent with Theorems~1--4. Compared with Algorithm~1, Algorithm~2 achieves consensus faster, exhibits greater stability and faster decrease in time-averaged dynamic regret, in this simulation. Fig.~3 further illustrates the trade-off between function-value sampling cost and convergence speed as the dimension $d$ increases.

As shown in Table~2, averaging the per-threshold reductions in the required rounds relative to the matched TS baselines yields $40.60\%$ and $39.55\%$ for Algorithm~1, and $57.97\%$ and $62.30\%$ for Algorithm~2, for regret and consensus, respectively, in Setting~I. In Setting~II, these reductions increase to $65.87\%$ and $60.71\%$ for Algorithm~1, and $85.29\%$ and $85.74\%$ for Algorithm~2; the consensus averages exclude NR entries. 
%Algorithm~2 requires fewer rounds than Algorithm~1 across all tested criteria. 
These results further demonstrate the effectiveness of the proposed mixing mechanism, especially with deterministic zeroth-order estimation, and
suggest larger relative gains under sparser trusted inter-cluster communication. Explaining this trend and extending the guarantees beyond uniform joint connectivity remain future work.

Fig.~1 also indicates that the two attack-free baselines fail to exhibit sublinear dynamic regret under Byzantine edge attacks. Moreover, Fig.~1 and Table~2 show that the trusted-skeleton algorithm reduces the dynamic regret more slowly than the proposed algorithms, providing empirical support for the usefulness of the proposed filtering mechanism in exploiting information from non-trusted edges. 
The adapted RDO-T algorithm, originally designed for node attacks, does not reach the prescribed regret threshold within the simulation horizon.
}

\section{Conclusion}
%In this paper, we proposed two types of zeroth-order resilient distributed online algorithms for networks under Byzantine edge attacks. Consensus of the estimates and sublinearity of the dynamic regret are established. %For future research, we are further improving the boundaries of the regret function and we will focus on the distributed games under Byzantine attacks in network. 
This paper studies zeroth-order distributed online optimization over networks subject to Byzantine edge attacks, where gradient information is unavailable and communication links may be adversarially corrupted. To address these challenges, two classes of resilient zeroth-order distributed online algorithms with attack-resilient mixing coefficient designs are developed. 
{The algorithms achieve asymptotic consensus and sublinear dynamic regret under the specified connectivity, minimizer-path variation, and parameter conditions, } thereby enabling reliable online tracking in adversarial environments. A comparative analysis of the two algorithms reveals an inherent trade-off between the computational cost of function-value sampling and the convergence speed. Numerical simulations further validate the theoretical results. 
Future work will focus on tightening the dynamic regret bounds, extending the guarantees beyond uniform joint connectivity, and exploring distributed games under attacks.
%This paper develops two zeroth-order distributed online algorithms with resilient mixing for networks under Byzantine edge attacks. {The algorithms achieve asymptotic consensus and sublinear dynamic regret under the specified connectivity, minimizer-path variation, and parameter conditions.}
%Their comparison highlights a trade-off between function-value sampling cost and convergence performance. Simulations support the theoretical results. 
%Future work includes extending the guarantees beyond uniform joint connectivity, and exploring distributed games under attacks.

\renewcommand\thedefinition{A\arabic{definition}}
\makeatother\setcounter{definition}{0}

\renewcommand\thelemma{A\arabic{lemma}}
\makeatother\setcounter{lemma}{0}

\appendix
\renewcommand{\thesection}{Appendix \Alph{section}}
\renewcommand{\thesubsection}{\Alph{section}.\arabic{subsection}}
\renewcommand{\thesubsubsection}{\Alph{section}.\arabic{subsection}.\arabic{subsubsection}}
\renewcommand{\theequation}{\Alph{section}.\arabic{equation}}
\renewcommand{\thefigure}{\Alph{section}.\arabic{figure}}
\renewcommand{\thetable}{\Alph{section}.\arabic{table}}

\section{Commonly used results}
%In what follows, we will present commonly used lemmas in graph theory and convex analysis. 

%~{\it{assumption A1:}}~For the adjacent matrices $A_t,~t\geqslant 0$,
%\begin{itemize}
%\item[a)] There exists a constant $ \eta > 0 $ such that for any $ i \in V $, if there is information transmission from node $ j $ to node $ i $ at time $ t $, then $ a^{ij}_t \geqslant \eta $; otherwise, $ a^{ij}_t = 0 $;

%\item[b)] There exists a positive integer $ B \geqslant 1 $ such that for any $ t \geqslant 0 $, $ (V,~E_t \cup E_{t+1} \cup \ldots \cup E_{t+B-1})$ is strongly connected.
%\end{itemize}

For a time-varying network graph $ \mathcal{G}_t = (V, E_t, A_t) $, define the transition matrix of the network by
$\Phi_{t,s}\triangleq A_{t}A_{t-1}\cdots A_{s+1}A_{s}\text{, }\forall t\geqslant s$ and $\Phi_{t,t+1}\triangleq I$. 

\begin{lemma}[see \cite{2009Nedic,2018Xie}] \label{Lem2.1} 
%若假设 \textcolor{black}{a) -- d)} 成立，那么有以下结论：
Assume $\{A_t,~t\geqslant 0\}$ is a row-stochastic matrix sequence. For $i\in V$ and $t\ge 0$, $(i,i)\in E_t$, there exist $\eta>0$,~positive integer $B \geqslant 1$ s.t. $ a^{ij}_t \geqslant \eta >0$ if $(i,j)\in E_t$ and $ a^{ij}_t = 0 $ otherwise, and $ (V,~E_t \cup E_{t+1} \cup \ldots \cup E_{t+B-1})$ is strongly connected. For the transition matrix $\Phi_{t,s},~t\geqslant s\geqslant 0$, there exists a sequence $\left\{\pi_t\in \mathbb{R}^n, {t} \geqslant 0\right\}$, such that
%$\left\{\pi_t=[\pi_t[1], \pi_t[2], \ldots, \pi_t[n]^{T}, {t} \geqslant 0\right\}$
\begin{itemize}
\item[a)] $1^{\top} \pi_t=1$ holds for $t\geqslant 0$, i.e., $\pi_t$ is  stochastic. 
\item[b)] There exist positive constants $C>0$ and $\lambda \in(0,1)$, such that for $t\geqslant s\geqslant 0$
\begin{align}
   \mid\Phi_{t,s}[i,j]-\pi_s[j]\mid \le C\lambda^{t-s}.\label{Phi}
\end{align}
%其中 $ \gamma=\frac{2(1+\eta^{-B} )}{1-\eta^{B} }$，$ \beta=(1-\eta^{B})^{1/B}$，$B=(m-1)B_0$，$ \eta $， $B_{0}$ 是上述假设中的参数。
\item[c)]There exists a positive constant $r_0>0$ such that $\pi_t[i] \geqslant r_0$ holds for  $i \in V$ and $t \geqslant 0$;
\item[d)] $\pi^{\top}_t=\pi^{\top}_{t+1} A_t$ holds for $t \geqslant 0$. 
\end{itemize}
\end{lemma}

%A set $ S \subset \mathbb{R}^d $ is called a convex set in $ \mathbb{R}^d $, if for any $ x, y \in S $ and  $ \lambda \in [0,1] $ it follows that $ \lambda x + (1 - \lambda) y \in S $. 

%\begin{definition}(see \cite{2012Shalevshwartz})%,1990Harker,2012Hazan
%\label{Def_conv} Denote by $\text{dom} (f)\subset \mathbb{R}^d$ the domain of $f(\cdot)$.
%A vector-valued function $\partial f(x)\in \mathbb{R}^d$ is called the subgradient of a non-smooth convex function $f(x):\text{dom} (f)\to \mathbb{R}$ if for any $x,y\in \text{dom}(f)$,
%$$f(y)-f(x)\geqslant \partial f(x)^\top (y-x).$$ 
%$f(\cdot)$ is called convex on $X\subset \text{dom} (f)$, if there exists a subgradient $\partial f(x)$ for any $x\in X$. 
%$f(\cdot)$ is called $\sigma$-strongly convex on $X\subset \text{dom} (f)$, if for any $x,y\in X$, $$f(y)-f(x)\geqslant \partial f(x)^\top (y-x)+\frac{\sigma}{2}\|y-x\|^2.$$
%\end{definition}

\begin{lemma}[see \cite{2010Ram,2022Xiong}] \label{Lem2.3} 
For the projection operator $\mathcal{P}_{X} (\cdot):\mathbb{R}^d \to X$, it holds 
\begin{align}
&\|\mathcal{P}_{X}(x)-\mathcal{P}_{X}(y)\|\le \|x-y\|,~\forall x,y \in \mathbb{R}^d \label{Proj1}
\\ & \left\langle\mathcal{P}_{X}(x)-x,\mathcal{P}_{X}(x)-y\right\rangle\le 0,~\forall x\in \mathbb{R}^d,y\in X.\label{Proj2}
\end{align}

%For $x^{1}, \ldots , x^{m} \in \mathbb{R}^d$,
%\begin{equation}\label{Avar}
%\sum_{i=1}^n\left\|x^{i}-\frac{1}{m}\sum_{j=1}^m x^{j}\right\|\le \sum_{i=1}^n\left\|x^{i}-x\right\|,~~\forall x\in \mathbb{R}^d.
%\end{equation}
\end{lemma}

\section{Proof of Lemma \ref{Lem:pi}}\label{Proof:Lem pi}
{By Assumption~3, \eqref{B-lower}, $B_t[i,j]=0$ for nontrusted $j$, as well as the push-sum scaling update \eqref{z}, Proposition~1 in \cite{2019Scutari} guarantees (\ref{theta}) in a). As to b), by {\eqref{A-row},} $\tilde{A}_{t}[i,j]\geqslant (\xi_1\xi_2)^{-1}{\beta_B}$ for $(i,j)\in {E}_t^r$ {and $\tilde{A}_{t}[i,j]=0$ otherwise}, together with Assumption 3, Lemma \ref{Lem2.1} yields b)}. \hfill$\square$

\section{Proof of Theorem \ref{sThm1}}\label{Proof:Thm 1.1}
By Assumption 2, (\ref{sa3}) and $u_t^i$ is uniformly distributed over $\mathbb{S}$, we have 
\begin{align}
\nonumber&\|g_t^i\|\le\frac{d}{2\delta_t}\left\|f_t^i(y_t^i+\delta_t u_t^i)-f_t^i(y_t^i-\delta_t u_t^i)\right\|\left\|u_t^i\right\|
\\
&\le \frac{d}{2\delta_t}\cdot L\cdot2\delta_t=dL\triangleq L_g. \label{s3a}
\end{align}
By Assumption 1, 4,  and $\textbf{0}\in X$, it's not hard to get 
%\begin{align*}
%&\left(1-\frac{\delta_{t-1}}{r}\right)x=\left(1-\frac{\delta_{t}}{r}\right)\cdot\frac{1-\frac{\delta_{t-1}}{r}}{1-\frac{\delta_{t}}{r}}x
%\\=&\left(1-\frac{\delta_{t}}{r}\right)\left(\frac{1-\frac{\delta_{t-1}}{r}}{1-\frac{\delta_{t}}{r}}x+\left(1-\frac{1-\frac{\delta_{t-1}}{r}}{1-\frac{\delta_{t}}{r}}\right)\textbf{0}\right)
%\\\in&\left(1-\frac{\delta_{t}}{r}\right)X=X_{\delta_{t}}, 
%\end{align*}
%with $\textbf{0}\in X\subset \mathbb{R}^d$, and thus for $t\geqslant 2$
\begin{align*}
X_{\delta_{t-1}}\subset X_{\delta_{t}},~t\geqslant 1,
\end{align*}
from which and by (\ref{A})- (\ref{A-row}), (\ref{sa1}) and (\ref{sa4}), we have 
\begin{align}
y_t^i=\sum_{j=1}^n \tilde{A}_t[i,j]x_t^j\in X_{\delta_{t-1}}\subset X_{\delta_{t}}.\label{s4a}
\end{align}
For $t\geqslant 0$, define $q^i_{t+1}=x^i_{t+1}-y^i_t$. By (\ref{sa1})-(\ref{sa4}), (\ref{s4a}), Lemma \ref{Lem2.3},  (\ref{theta}) and (\ref{s3a}), it follows that 
\begin{align}
\nonumber\left\|q^i_{t+1}\right\|&=  \left\|\mathcal{P}_{X_{\delta_{t}}}\left(y^i_t-\alpha_t \frac{g_t^i}{ z_{t+1}[i] }\right)-\mathcal{P}_{X_{\delta_{t}}}\left(y^i_t\right)\right\| \\
& \le \frac{\alpha_t}{ z_{t+1}[i] }\left\|g_t^i\right\| \le \xi_1 L_g \alpha_t\label{s4}
\end{align}
By (\ref{sa1}) and the definition of $q^i_{t+1}$, for $t \geqslant 2$ we have 
\begin{align}
x^i_{t+1} %=y^i_t+q^i_{t+1}
\!=\!\sum_{j=1}^n \!\tilde{\Phi}_{t, 1}[i, j] {x}^j_1\!+\!\sum_{s=2}^t \sum_{j=1}^n \!\tilde{\Phi}_{t, s}[i, \!j] {q}^{j}_s\!+\!q^i_{t+1}, \label{s6}
\end{align}
%On the other hand, 
and the definition of $\bar{x}_{t+1}$,  $q^i_{t+1}$, $y_t^i$ and (\ref{Phi'2}) yields 
\begin{align}
\nonumber& \bar{x}_{t\!+\!1} \!\triangleq \!\sum_{i=1}^n \!\pi_{t\!+\!1}[i] x^i_{t\!+\!1}\!=\!\sum_{i=1}^n \!\pi_{t\!+\!1}[i] y_t^i\!+\!\sum_{i=1}^n\! \pi_{t\!+\!1}[i] q^i_{t\!+\!1}\nonumber
\\
%& =\sum_{i=1}^n \pi_{t+1}[i] \sum_{j=1}^n \tilde{A}_t[i,j] x^j_{t}+\sum_{i=1}^n \pi_{t+1}[i] q^i_{t+1} \nonumber\\
& =\sum_{j=1}^n\left(\sum_{i=1}^n \pi_{t+1}[i] \tilde{A}_t[i,j]\right) x^j_{t}+\sum_{i=1}^n \pi_{t+1}[i] q^i_{t +1}\nonumber\\
& =\sum_{j=1}^n \!\pi_{t}[j] x^j_{t}\!\!+\!\!\sum_{i=1}^n\!\pi_{t\!+\!1}[i] q^i_{t\!+\!1}\!\!=\!\!\bar{x}_{t}\!\!+\!\!\sum_{i=1}^n \!\pi_{t\!+\!1}[i] q^i_{t\!+\!1}, \label{s7}
\end{align}
from which, by induction we have
\begin{align}
\bar{x}_{t+1} =\sum_{j=1}^n \pi_1[j] x^j_1+\sum_{s=2}^{t+1} \sum_{j=1}^n \pi_s[j] q^j_s\label{s7b}
\end{align}
Combining \eqref{s6} and \eqref{s7b}, together with (\ref{Phi'1}), (\ref{s4}), Assumption 1, and $\pi_t$ is stochastic,  we can obtain 
\begin{align*}
& \nonumber\left\|x^i_{t+1}-\overline{x}_{t+1}\right\|
 \nonumber\le \sum_{j=1}^{n}\left\|\tilde{\Phi}_{t, 1}[i, j]-\pi_1[j]\right\| \cdot\left\|x^j_1\right\| \\
\nonumber+&\sum_{s=2}^t \!\sum_{j=1}^{n}\!\left\|\tilde{\Phi}_{t,s}[i, \!j]\!\!-\!\!\pi_s[j]\right\| \left\|{q}^j_s\right\|\!+\!\sum_{j=1}^n \!\pi_{t\!+1\!}[j]\|q^j_{t\!+\!1}\|\!\!+\!\!\|q^i_{t\!+\!1}\|\\
\le & {nCM} \lambda^{t\!-\!1}\!+\!{n}\xi_1 L_g {C}\sum_{s=2}^{t} \lambda^{t\!-\!s} {\alpha_{s\!-\!1}}\!+\!2\xi_1L_g\alpha_t  \!\triangleq\! \beta_{t+1},
\end{align*}
which yields Theorem \ref{sThm1} i). If $\lim_{t\rightarrow \infty}\alpha _{t}=0$, then by $\lambda\in(0,1)$, it is direct to prove that $\lim_{t\rightarrow \infty}{\beta}_{t}=0$.\hfill$\square$

\section{Proof of Lemma \ref{sLem4.2}}\label{Proof:Lem 1.1}
{Define $\tilde{x}_t^*\triangleq(1-\delta_t/r)x_t^*\in X_{\delta_t}$ and $\theta_t=\left\|x_t^*-x_{t+1}^*\right\|$. 

By Assumption \ref{asmp:X-interior}, \ref{asmp:delta},} it can be easily verified that 
{\begin{align}
\nonumber&~~ \frac{1}{2}\left\|\bar{x}_{t+1}-\tilde{x}_{t+1}^*\right\|^2-\frac{1}{2}\left\|\bar{x}_t-\tilde{x}_t^*\right\|^2 
%\\\nonumber= & \frac{1}{2}\left\langle\bar{x}_{t+1}+\bar{x}_t-{x}_{t+1}^*-x_t^*, \bar{x}_{t+1}-\bar{x}_t-x_{t+1}^*+x_t^*\right\rangle \\
%\nonumber= & \frac{1}{2}\left\langle 2 \bar{x}_{t+1}\!-\!x_{t+1}^*\!-\!x_t^*\!+\!\bar{x}_t\!-\!\bar{x}_{t+1}, \bar{x}_{t+1}-\bar{x}_t-x_{t+1}^*+x_t^*\right\rangle \\
%= & -\frac{1}{2}\left\|\bar{x}_t-\bar{x}_{t+1}\right\|^2+\frac{1}{2}\left\langle\bar{x}_t-\bar{x}_{t+1}, x_t^*-x_{t+1}^*\right\rangle \\
%\nonumber& +\frac{1}{2}\left\langle\bar{x}_t-\bar{x}_{t+1}, x_{t+1}^*+x_t^*-2 \bar{x}_{t+1}\right\rangle
%\\\nonumber&+\frac{1}{2}\left\langle 2 \bar{x}_{t+1}-x_{t+1}^*-x_t^*, x_t^*-x_{t+1}^*\right\rangle 
\\
\nonumber&=-\frac{1}{2}\left\|\bar{x}_t-\bar{x}_{t+1}\right\|^2+\left\langle\bar{x}_t-\bar{x}_{t+1}, \tilde{x}_t^*-\bar{x}_{t+1}\right\rangle
\\\nonumber&+\frac{1}{2}\left\langle 2 \bar{x}_{t\!+\!1}\!-\!\tilde{x}_{t+1}^*\!-\!\tilde{x}_t^*, \tilde{x}_t^*\!-\!\tilde{x}_{t\!+\!1}^*\right\rangle\!\le\! 2M\theta_t\! -\!\frac{1}{2}\left\|\bar{x}_t\!-\!\bar{x}_{t\!+\!1}\right\|^2\!\\&+\left\langle\bar{x}_t\!-\!\bar{x}_{t\!+\!1}, \tilde{x}_t^*\!-\!\bar{x}_{t\!+\!1}\right\rangle\!+\!\frac{2M^2}{r}(\delta_t-\delta_{t+1}), \label{s9}
\end{align}}
{where the last inequality uses
\begin{align*}
\|\tilde{x}_{t+1}^*-\tilde{x}_t^*\|&=\left\|\left(1\!-\!\frac{\delta_{t+1}}{r}\right)({x}_{t+1}^*-{x}_t^*)+\frac{\delta_t\!-\!\delta_{t+1}}{r}{x}_t^*\right\|
\\&\leq \theta_t+\frac{M}{r}(\delta_t-\delta_{t+1})
\end{align*}}
{By (\ref{s7}), (\ref{s4}), Schwarz inequality, (\ref{theta}), Theorem \ref{sThm1} and $\pi_t$ is stochastic, we have 
\begin{align}
\nonumber& \left\langle\bar{x}_{t+1}-\bar{x}_t, \bar{x}_{t+1}-{\tilde{x}_t^*}\right\rangle 
%\\\nonumber =&\left\langle\sum_{i=1}^n \pi_{t+1}[i] q_{t+1}^i, \bar{x}_{t+1}-x_t^*\right\rangle 
%\end{align}\begin{align}
\\
\nonumber =&\sum_{i=1}^n \!\!\pi_{t\!+\!1}[i]\!\left\langle q_{t\!+\!1}^i, \bar{x}_{t\!+\!1}\!\!-\!\!x_{t\!+\!1}^i\right\rangle\!\!+\!\!\sum_{i=1}^n\! \pi_{t\!+\!1}[i]\!\left\langle q_{t\!+\!1}^i, x_{t\!+\!1}^i\!-\!\tilde{x}_t^*\right\rangle \\
 \le& \xi_1 L_g \alpha_t \beta_{t\!+\!1}\!+\!\sum_{i=1}^n \pi_{t\!+\!1}[i]\left\langle q_{t\!+\!1}^i, x_{t+1}^i\!-\!\tilde{x}_t^*\right\rangle\label{s11a}
\end{align}
Next we focus on $\pi_{t+1}[i]\left\langle q_{t+1}^i, x_{t+1}^i-\tilde{x}_t^*\right\rangle$. Note that 
\begin{align}
&\nonumber\pi_{t\!+\!1}[i]\left\langle\! q_{t\!+\!1}^i, x_{t\!+\!1}^i\!\!-\!\!\tilde{x}_t^*\right\rangle
\!\!=\!\!\pi_{t\!+\!1}[i]\left\langle q_{t\!+\!1}^i\!+\!\alpha_t\frac{g_t^i}{z_{t\!+\!1}[i]}, x_{t\!+\!1}^i\!-\!\tilde{x}_t^*\!\right\rangle\\&+
\pi_{t+1}[i]\left\langle \alpha_t\frac{g_t^i}{z_{t+1}[i]}, \tilde{x}_t^*-x_{t+1}^i\right\rangle.
\end{align}
On the one hand, (\ref{Proj2}) in Lemma \ref{Lem2.3} yields 
\begin{align}
\nonumber&\pi_{t+1}[i]\left\langle q_{t+1}^i+\alpha_t\frac{g_t^i}{z_{t+1}[i]}, x_{t+1}^i-\tilde{x}_t^*\right\rangle
\\=&\pi_{t\!+\!1}[i]\left\langle\! x_{t\!+\!1}^i\!\!-\!\!\left(y_t^i\!\!-\!\!\alpha_t\frac{g_t^i}{z_{t\!+\!1}[i]}\right), x_{t\!+\!1}^i\!-\!\tilde{x}_t^*\!\right\rangle\!\le\! 0.\label{s11}
\end{align}
On the other hand, by Cauchy inequality, Theorem \ref{sThm1}, (\ref{theta}), (\ref{s3a}) and $\pi_t$ is stochastic, it follows that
\begin{align}
\nonumber&\sum_{i=1}^n\pi_{t+1}[i]\left\langle \alpha_t\frac{g_t^i}{z_{t+1}[i]}, \tilde{x}_t^*-x_{t+1}^i\right\rangle
\\\nonumber=&\sum_{i=1}^n\pi_{t\!+\!1}[i]\left\langle \alpha_t\frac{g_t^i}{z_{t\!+\!1}[i]}, \tilde{x}_t^*\!-\!\bar{x}_t\!+\! \bar{x}_t\!-\!\bar{x}_{t\!+\!1}\!+\!\bar{x}_{t\!+\!1}\!-\!x_{t\!+\!1}^i\right\rangle
\\\nonumber\leq&\sum_{i=1}^n\pi_{t+1}[i]\left\langle \alpha_t\frac{g_t^i}{z_{t+1}[i]}, \tilde{x}_t^*-\bar{x}_t\right\rangle
\\+& \frac{1}{2}\xi_1^2\alpha_t^2L_g^2\!+\!\frac{1}{2}\left\|\bar{x}_t\!-\!\bar{x}_{t\!+\!1}\right\|^2\!+\!\xi_1 L_g\alpha_t\beta_{t\!+\!1}\label{s12}
\end{align}
By $\eqref{row}$, we can get $z_s^\top=z_{t+1}^{\top}\widetilde\Phi_{t,s}$ for any fixed $s\geq0$ and $t\geq s$. By $z_t$ is stochastic and Lemma \ref{Lem:pi}, we have
\begin{align*}
|z_s[j]-\pi_s[j]|
\le \sum_i z_{t+1}[i]\,
|\widetilde\Phi_{t,s}[i,j]-\pi_s[j]|
\le C\lambda^{t-s}, 
\end{align*}
Letting $t\to\infty$ gives
\begin{align}
\pi_s=z_s,~\forall s\geq0,\label{eq:pi=s}
\end{align}
Combining (\ref{s11a}) -- (\ref{eq:pi=s}), we can get
\begin{align}
\nonumber&\left\langle\bar{x}_{t+1}\!-\!\bar{x}_t, \bar{x}_{t+1}\!-\!\tilde{x}_t^*\right\rangle
\!\le \!\sum_{i=1}^n\left\langle \alpha_t{g_t^i}, \tilde{x}_t^*\!-\!\bar{x}_t\right\rangle
\\&+ \frac{1}{2}\xi_1^2\alpha_t^2L_g^2
+\frac{1}{2}\left\|\bar{x}_t\!-\!\bar{x}_{t\!+\!1}\right\|^2\!+\!2\xi_1 L_g\alpha_t\beta_{t\!+\!1}, 
\end{align}
which together with (\ref{s9}) yields 
\begin{align}
&\nonumber\frac{1}{2}\left\|\bar{x}_{t+1}\!-\!\tilde{x}_{t+1}^*\right\|^2\!-\!\frac{1}{2}\left\|\bar{x}_t\!-\!\tilde{x}_t^*\right\|^2
\le 2M\theta_t\!+\!2\xi_1L_g\alpha_t\beta_{t\!+\!1}
\\+&\!\frac{1}{2}\xi_1^2\alpha_t^2L_g^2\!+\!\alpha_t\sum_{i=1}^n\!\left\langle g_t^i, \tilde{x}_t^*\!-\!\bar{x}_t\right\rangle\!+\!\frac{2M^2}{r}(\delta_t\!-\!\delta_{t\!+\!1})
\label{s13}
\end{align}
%Denote $\mathcal{F}_t=\sigma\left\{u_s^i:~0\le s\le t,~i\in V \right\}$, by (\ref{sa1}) and (\ref{sa4}) we have $x_t^i,~y_t^i \in \mathcal{F}_{t-1}$ for $t\geqslant 1$ and $i\in V$. Note that by (\ref{aa1'}) we can obtain
%\begin{align}
%\mathbb{E}\left(g_t^i|y_t^i\right)=\mathbb{E}\left(g_t^i|\mathcal{F}_{t-1}\right)=\nabla \hat{f}_{t,\delta_t}^i(y_t^i),
%\end{align}
Next we focus on $\sum_{i=1}^n\left\langle g_t^i, \tilde{x}_t^*-\bar{x}_t\right\rangle$. By (\ref{aa2'}), (\ref{sa3}), and the fact that $u_t^i$ independent of $\sigma(y_t^i)$,  we can obtain
\begin{align}
\mathbb{E}\left(g_t^i|y_t^i\right)=\nabla\hat{f}_{t,\delta_t}^i(y_t^i). \label{eq:add-ex}
\end{align}
%which combines (\ref{sbar_x}) and $\bar{x}_t\in \mathcal{F}_{t-1}$  yields 
%\begin{align}
%\nonumber& \sum_{i=1}^n \frac{\pi_{t+1}[i]}{z_{t+1}[i]}\mathbb{E}\left\langle g_t^i, x_t^*-\bar{x}_t\right\rangle \\
%\nonumber= & \sum_{i=1}^n \frac{\pi_{t+1}[i]}{z_{t+1}[i]}\mathbb{E}\left[\mathbb{E}\left\langle g_t^i, x_t^*-\bar{x}_t\right\rangle|\mathcal{F}_{t-1}\right]\\
%\nonumber= & \sum_{i=1}^n \frac{\pi_{t+1}[i]}{z_{t+1}[i]}\mathbb{E}\left\langle \mathbb{E}\left(g_t^i|\mathcal{F}_{t-1}\right), x_t^*-\bar{x}_t\right\rangle\\
%\nonumber= & \sum_{i=1}^n \frac{\pi_{t+1}[i]}{z_{t+1}[i]}\mathbb{E}\left\langle\nabla \hat{f}_{t,\delta_t}^i\left(y_t^i\right), x_t^*-\bar{x}_t\right\rangle \\
%\nonumber= & \sum_{i=1}^n \frac{\pi_{t+1}[i]}{z_{t+1}[i]}\mathbb{E}\left\langle\nabla \hat{f}_{t,\delta_t}^i\left(y_t^i\right), x_t^*-y_t^i\right\rangle+\sum_{i=1}^n \frac{\pi_{t+1}[i]}{z_{t+1}[i]}\mathbb{E}\left\langle\nabla \hat{f}_{t,\delta_t}^i\left(y_t^i\right), y_t^i-\bar{x}_t\right\rangle\\=&E(1)+E(2), \label{s14}
%\end{align}
Since $\tilde{x}_t^*$ is deterministic and $\tilde{x}_t^*-y_t^i$ is $\sigma(y_t^i)-$measurable, \eqref{eq:add-ex} yields 
\begin{align}
\nonumber& \sum_{i=1}^n \mathbb{E}\left\langle g_t^i, \tilde{x}_t^*\!-\!\bar{x}_t\right\rangle =\sum_{i=1}^n\left(\mathbb{E}\left\langle g_t^i, \tilde{x}_t^*\!-\!y_t^i\!+\!y_t^i\!-\!\bar{x}_t\right\rangle\right)\\
\nonumber= & \sum_{i=1}^n \left(\mathbb{E}\!\left[\mathbb{E}\left\langle g_t^i, \tilde{x}_t^*-y_t^i\right\rangle|y_t^i\right]+\mathbb{E}\left\langle g_t^i, y_t^i-\bar{x}_t\right\rangle\right)\\
\nonumber= & \sum_{i=1}^n \left(\mathbb{E}\left\langle \mathbb{E}\left(g_t^i|y_t^i\right), \tilde{x}_t^*-y_t^i\right\rangle+\mathbb{E}\left\langle g_t^i, y_t^i-\bar{x}_t\right\rangle\right)
\\
\nonumber= & \sum_{i=1}^n \left(\mathbb{E}\left\langle\nabla \hat{f}_{t,\delta_t}^i\left(y_t^i\right), \tilde{x}_t^*-y_t^i\right\rangle\!+\mathbb{E}\left\langle g_t^i,  y_t^i-\bar{x}_t\right\rangle\right)
\\
=&E_t(1)+E_t(2), \label{s14}
\end{align}
with }
\begin{align*}
{E_t(1)=}& {\sum_{i=1}^n \mathbb{E}\left\langle\nabla \hat{f}_{t,\delta_t}^i\!\left(y_t^i\right), \tilde{x}_t^*-y_t^i\right\rangle,}\\
E_t(2)=& \sum_{i=1}^n \mathbb{E}\left\langle g_t^i, y_t^i-\bar{x}_t\right\rangle.  
\end{align*}
By (\ref{A}), (\ref{A-row}), (\ref{sa1}) and Theorem \ref{sThm1}, it follows that
\begin{align}
\nonumber&\left\|y_t^i\!-\!\bar{x}_t\right\|=\left\|\sum_{j=1}^n \!\tilde{A}_t[i,j]x_t^j\!-\!\bar{x}_t\right\|\!=\!\left\|\sum_{j=1}^n \!\tilde{A}_t[i,j]\left(x_t^j\!-\!\bar{x}_t\right)\right\|
\\&\leq\sum_{j=1}^n \tilde{A}_t[i,j]\left\|x_t^j-\bar{x}_t\right\|, \text{thus }\mathbb{E}\left\|y_t^i\!-\!\bar{x}_t\right\| \!\leq\!\beta_t, \label{s13b}
\end{align}
from which and by Schwarz inequality, (\ref{s3a}), for $E_t(2)$ we have 
\begin{align}
E_t(2)\le \sum_{i=1}^n \mathbb{E}\left[\left\|g_t^i\right\|\cdot \left\|y_t^i-\bar{x}_t\right\|\right]
\le nL_g\beta_t\label{s13bb}
\end{align}
Next we consider $E_t(1)$. Noting the definition of $\delta-$smoothing function, Assumption 1,2, and the Lebesgue’s dominated convergence theorem, it can be proven  $\nabla \hat{f}_{t,\delta_t}^i\!\left(x\right)\!=\!\mathbb{E}\left[\nabla f_t^i(x+\delta_t v)\right]$ for $x\in X_{\delta_t}$ and $\hat{f}_{t,\delta_t}^i\!\left(x\right)$ is also (strongly) convex over $X_{\delta_t}$ if $f_t^i(x)$ is (strongly) convex over $X$, 
%\mathbb{E}\left(\right)
thus for $x\in X_{\delta_t}$, by Assumption \ref{asmp:convex-feasible-set} we can get 
\begin{align}
\left\|\nabla\hat{f}_{t,\delta_t}^i\!\left(x\right)\right\|\!
%=&\left\|\mathbb{E}\left[\nabla f_t^i(x\!+\!\delta_t u)\right]\right\|\\
\le \mathbb{E}\left\|\nabla f_t^i(x+\delta_t v)\right\|\le L.  \label{s13ab}
\end{align}
By $\hat{f}_{t,\delta_t}^i(\cdot)$ is $\sigma-$strongly convex function on $X_{\delta_t}$, 
{\eqref{s13b}-\eqref{s13ab} and  Assumption \ref{asmp:X-interior}, we have
\begin{align}
\nonumber&\sum_{i=1}^n \left\langle\nabla \hat{f}_{t,\delta_t}^i\left(y_t^i\right), \tilde{x}_t^*-y_t^i\right\rangle
\\
\nonumber\le&\sum_{i=1}^n \left(\hat{f}_{t,\delta_t}^i\left(\tilde{x}_t^*\right)-\hat{f}_{t,\delta_t}^i\left(y_t^i\right)-\frac{\sigma}{2}\left\|y_t^i-\tilde{x}_t^*\right\|^2\right)
\\\nonumber\le &\sum_{i=1}^n \left(\hat{f}_{t,\delta_t}^i\left(\tilde{x}_t^*\right)-\hat{f}_{t,\delta_t}^i\left(\bar{x}_t\right)+\hat{f}_{t,\delta_t}^i\left(\bar{x}_t\right)-\hat{f}_{t,\delta_t}^i\left(y_t^i\right)\right.
\\
\nonumber&\left.-\frac{\sigma}{2}\left\|y_t^i-\tilde{x}_t^*\right\|^2\right)
\le \sum_{i=1}^n \left(\hat{f}_{t,\delta_t}^i\left(\tilde{x}_t^*\right)-\hat{f}_{t,\delta_t}^i\left(\bar{x}_t\right)\right.\\
&~\left.+L\beta_t-\frac{\sigma}{2}\left\|y_t^i-\tilde{x}_t^*\right\|^2\right)\label{s15bbb}, 
\end{align}
Schwarz inequality, Assumption 1 and (\ref{s13b}) yields 
\begin{align}
\nonumber&\mathbb{E}\left|\left\|y_t^i-\tilde{x}_t^*\right\|^2-\left\|\bar{x}_t-\tilde{x}_t^*\right\|^2\right|
\\=&\mathbb{E}\left|\left\langle y_t^i+\bar{x}_t-2\tilde{x}_t^*, y_t^i-\bar{x}_t\right\rangle\right|
%\le&\left(\|y_t^i\|+\|\bar{x}_t\|+2\|\tilde{x}_t^*\|\right)\|y_t^i-\bar{x}_t\| 
\le4M\beta_t. \label{s14bb}
\end{align}
%For $E_t(1)$, (\ref{s15})-(\ref{s15bbb}) yields 
%\begin{align}
%\nonumber&E_t(1)=\sum_{i=1}^n \frac{\pi_{t+1}[i]}{z_{t+1}[i]}\mathbb{E}\left\langle\nabla \hat{f}_{t,\delta_t}^i\left(y^i_t\right), x_t^*-y^i_t\right\rangle
%\\\nonumber \leq& \sum_{i=1}^n\mathbb{E}\left[\hat{f}_{t,\delta_t}^i\left(x_t^{*}\right)-\hat{f}_{t,\delta_t}^i\left(\bar{x}_t\right)-\frac{\sigma}{2}\left\|\bar{x}_t-x_t^*\right\|^2\right]
%\\&+2MnL\xi_1C\lambda^{T-t-2}+(4M+L)n\beta_t. \label{s22}
%\end{align}
Combining (\ref{s13bb})-(\ref{s14bb}), for (\ref{s14}) we have
\begin{align}
\nonumber&\sum_{i=1}^n \mathbb{E}\left\langle g_t^i, \tilde{x}_t^*-\bar{x}_t\right\rangle
\leq\sum_{i=1}^n\mathbb{E}\left[\hat{f}_{t,\delta_t}^i\left(\tilde{x}_t^*\right)-\hat{f}_{t,\delta_t}^i\left(\bar{x}_t\right)\right.
\\&~\left.-\frac{\sigma}{2}\left\|\bar{x}_t-\tilde{x}_t^*\right\|^2\right]+\left(\left(2M\sigma+L\right)n+nL_g\right)\beta_t. 
\label{s17}
\end{align}}
Note that by the definition of $\delta-$smoothing function and Assumption 2, for $x\in X_{\delta_t}$ we can get 
\begin{align}
\left|\hat{f}_{t,\delta_t}^i\!(x)\!\!-\!\!f_t^i(x)\right|\!\!=\!\!\left|\mathbb{E}[f^i_t(x\!\!+\!\!\delta_t v)] \!\!-\!\!\mathbb{E}[f_t^i(x)]\right|\le L\delta_t\label{s17b}
\end{align}
where $v$ is a random vector uniformly distributed over $\mathbb{B}$. Combining (\ref{s17}), (\ref{s17b}) and  
(\ref{s13}), {and note that $|f_t^i(x_t^*)-f_t^i(\tilde{x}_t^*)|\leq ML\delta_t/r$}, we can obtain (\ref{s8}). \hfill$\square$

\section{Proof of Theorem \ref{sThm2}} \label{Proof:Thm 1.2}
If $f_t^i(\cdot)$ is $\sigma-$strongly convex on $X$, Lemma \ref{sLem4.2} yields{
\begin{align}
\nonumber&\sum_{t=1}^T \sum_{i=1}^n\mathbb{E}\left[f_t^i\left(\overline{x}_t\right)-f_t^i\left(x_t^*\right)\right]\\
\nonumber\le &\sum_{t=1}^T\!\mathbb{E}\!\left(\!\frac{\left\|\bar{x}_t\!\!-\!\!\tilde{x}_t^*\right\|^2\!\!-\!\!\left\|\bar{x}_{t\!+\!1}\!\!-\!\!\tilde{x}_{t\!+\!1}^*\right\|^2}{2\alpha_t}\!\right)\!\!-\!\!\frac{\sigma n}{2} \!\sum_{t=1}^T \!\mathbb{E}\left\|\bar{x}_t\!\!-\!\!\tilde{x}_t^*\right\|^2\\
\nonumber&+\sum_{t=1}^T\left(2 M \frac{\theta_t}{\alpha_t}+2\xi_1 L_g \beta_{t+1}+\frac{1}{2} \xi_1^2 L_g^2 \alpha_t+{C_4} \delta_t\right.
\\&\left.{+\frac{2M^2}{r}\frac{\delta_t-\delta_{t+1}}{\alpha_t}}+n\left(2M\sigma\!+\!L\!+\!L_g\right)\beta_t\right), \label{s23}
\end{align}
where $C_4\triangleq n\!\left(2L\!+\!LM/r\right)$.} Note that 
\begin{align}
\nonumber& \sum_{t=1}^T \frac{\left\|\bar{x}_t-\tilde{x}_t^*\right\|^2}{2\alpha_t}-\sum_{t=1}^T \frac{\left\|\bar{x}_{t+1}-\tilde{x}_{t+1}^*\right\|^2}{2\alpha_t} 
%\\\nonumber= & \frac{\left\|\bar{x}_{1}\!-\!x_{1}^*\right\|^2}{2\alpha_1}\!\!+\!\!\frac{1}{2}\!\sum_{t=2}^T \left\|\bar{x}_{t}\!-\!x_{t}^*\right\|^2\!\left(\frac{1}{\alpha_t}\!-\!\frac{1}{\alpha_{t\!-\!1}}\right)\!\!-\!\!\frac{\left\|\bar{x}_{T\!+\!1}\!-\!x_{T\!+\!1}^*\right\|^2}{2\alpha_T} 
\\
\le& \frac{\left\|\bar{x}_{1}\!-\!\tilde{x}_{1}^*\right\|^2}{2\alpha_1}\!+\!\frac{1}{2}\sum_{t=2}^T \left\|\bar{x}_{t}\!-\!\tilde{x}_{t}^*\right\|^2\left(\frac{1}{\alpha_t}\!-\!\frac{1}{\alpha_{t\!-\!1}}\right), \label{s24}
\end{align}
which combines (\ref{s23})  and Assumption 1 yields
{\begin{align}
\nonumber& \sum_{t=1}^T \sum_{i=1}^n\mathbb{E}\left[f_t^i\left(\overline{x}_t\right)-f_t^i\left(x_t^*\right)\right]
%\le \frac{\mathbb{E}\left\|\bar{x}_{1}-x_{1}^*\right\|^2}{2\alpha_1}
%\\\nonumber&+\!\frac{1}{2}\!\sum_{t=2}^T \mathbb{E}\left\|\bar{x}_{t}\!-\!x_{t}^*\right\|^2\left(\frac{1}{\alpha_t}\!-\!\frac{1}{\alpha_{t\!-\!1}}\right)\!-\!\frac{\sigma n}{2} \!\sum_{t=2}^T \mathbb{E}\left\|\bar{x}_t\!-\!x_t^*\right\|^2\\
%\nonumber& +\sum_{t=1}^T\left(2 M \frac{\theta_t}{\alpha_t}+2\xi_1 L_g \beta_{t+1}+\frac{1}{2} \xi_1^2 L_g^2 \alpha_t\right.
%\\\nonumber&\left.+\left(\left(4M\!+\!L\right)n\!+\!L_g\xi_1\right)\beta_t\!+\!2nL\delta_t\!+\!2MnL\xi_1C\lambda^{T-t-2}\right)
\\
\nonumber\le & \frac{2M^2}{\alpha_1}+\frac{1}{2}\sum_{t=2}^T \mathbb{E}\left\|\bar{x}_t-\tilde{x}_t^*\right\|^2\left(\frac{1}{\alpha_t}-\frac{1}{\alpha_{t-1}}-\sigma n\right)
\\
\nonumber&+\sum_{t=1}^T\left(2 M \frac{\theta_t}{\alpha_t}+2\xi_1 L_g \beta_{t+1}+\frac{1}{2} \xi_1^2 L_g^2 \alpha_t+\!{C_4}\delta_t\right.
\\&\left.{+\frac{2M^2}{r}\frac{\delta_t-\delta_{t+1}}{\alpha_t}}+n\left(2M\sigma\!+\!L\!+\!L_g\right)\beta_t\right).\label{s25}
\end{align}}
Since $f_t^i(\cdot)$ is convex on $X$, by Assumption 2, Theorem \ref{sThm1} and (\ref{s25}), we can get 
\begin{align}
&\mathbb{E}\left[\mathrm{regret}_T^{D,j}\right]
=  \sum_{t=1}^T \nonumber\sum_{i=1}^n\mathbb{E}\left[f_t^i\left(x_t^j\right)-f_t^i\left(x_t^*\right)\right] 
\\
\nonumber= & \!\sum_{t=1}^T\!\sum_{i=1}^n\!\mathbb{E}\!\left[f_t^i(x_t^j)\!\!-\!\!f_t^i\left(\bar{x}_t\right)\right]\!\!+\!\!\sum_{t=1}^T \!\sum_{i=1}^n\!\mathbb{E}\!\left[f_t^i\left(\bar{x}_t\right)\!\!-\!\!f_t^i\left(x_t^*\right)\right]\\
\le & n L \sum_{t=1}^T \beta_t+\sum_{t=1}^T \sum_{i=1}^n\mathbb{E}\left[f_t^i\left(\bar{x}_t\right)-f_t^i\left(x_t^*\right)\right] 
%\\\nonumber\le & \frac{2M^2}{\alpha_1}+\frac{1}{2}\sum_{t=2}^T \mathbb{E}\left\|\bar{x}_t-x_t^*\right\|^2\left(\frac{1}{\alpha_t}-\frac{1}{\alpha_{t-1}}-\sigma n\right) \\
%\nonumber& + \sum_{t=1}^T\left(2 M \frac{\theta_t}{\alpha_t}+\frac{1}{2} \xi_1^2 L_g^2 \alpha_t+2nL\delta_t\right) \\
%& +\left(4Mn+2n L+3L_g\xi_1\right) \sum_{t=1}^{T+1} \beta_t+\frac{2MnL\xi_1C}{\lambda^2(1-\lambda)}
\label{s26}
\end{align}
By (\ref{s86}) and Assumption 1, it holds that
\begin{align}
\nonumber&\sum_{t=3}^{T+1} \beta_t %=nCM\!\sum_{t=3}^{T+1} \!\lambda^{t\!-\!2}\!+\!\xi_1^2 C L_g \sum_{t=3}^{T+1} \sum_{s=2}^{t-1} \alpha_s \lambda^{t\!-\!s\!-\!1}\!\\
%\nonumber&~~+2\xi_1 L_g \sum_{t=3}^{T+1} \alpha_{t-1} \\
\!\le n\xi_1 C L_g\sum_{s=2}^T {\alpha_{s\!-\!1}}\!\! \!\sum_{t=s+1}^{\infty}\! \lambda^{t\!-\!s\!-\!1}\!\!+\!\!2\xi_1 L_g \sum_{t=3}^{T\!+\!1} \alpha_{t\!-\!1}\\
&~+\!\frac{n C M}{1\!-\!\lambda} \!\le\!\frac{n C M}{1-\lambda}\!+\!\left(\frac{n\xi_1 C L_g}{1\!-\!\lambda}+2\xi_1 L_g\right) \sum_{t=1}^T \alpha_t. \label{s27}
\end{align}
Combining (\ref{s25}) -- (\ref{s27}) we can get (\ref{s18}). (\ref{s18aa}) and (\ref{s18b}) can be directly obtained by (\ref{s18a}) and (\ref{s18}), respectively. 

Note that $\{\alpha_t\}_{t\geqslant 1}$ is a positive decreasing sequence yields $\left(\frac{1}{\alpha_t}-\frac{1}{\alpha_{t-1}}\right)\geqslant 0$. Taking $\sigma=0$ in (\ref{s25})-(\ref{s26}) yields
{\begin{align}
\nonumber&\mathbb{E}\left[\mathrm{regret}_T^{D,j}\right]\le \frac{2M^2}{\alpha_1}+2M^2\sum_{t=2}^T \left(\frac{1}{\alpha_t}-\frac{1}{\alpha_{t-1}}\right) 
\\
\nonumber&+\sum_{t=1}^T\left(2 M \frac{\theta_t}{\alpha_t}+\frac{1}{2} \xi_1^2 L_g^2 \alpha_t+{C_4}\delta_t{+\frac{2M^2}{r}\frac{\delta_t\!-\!\delta_{t\!+\!1}}{\alpha_t}}\right) \\
\nonumber& +\left({2\xi_1 L_g}+nL_g+2n L\right) \sum_{t=1}^{T+1} \beta_t
\\
\nonumber= & \frac{2M^2}{\alpha_T}\!+\! \sum_{t=1}^T\left(2 M \frac{\theta_t}{\alpha_t}\!+\!\frac{1}{2} \xi_1^2 L_g^2 \alpha_t\!+\!{C_4}\delta_t{+\frac{2M^2}{r}\frac{\delta_t\!-\!\delta_{t\!+\!1}}{\alpha_t}}\right) 
\\
&+ \left({2\xi_1 L_g}+nL_g+2n L\right) \sum_{t=1}^{T+1} \beta_t, \label{s29}
\end{align}}
which combines (\ref{s27}) yields (\ref{s18a}). \hfill$\square$

\section{Proof of Theorem \ref{Thm1}--\ref{Thm2}}\label{Proof:Thm2}
The proof is similar to that of Theorem \ref{sThm1}, Lemma \ref{sLem4.2} and Theorem \ref{sThm2}. Here we only list some key points. 

By Taylor's expansion and noting that $f_t^i(\cdot)$ is second order differentiable, it follows that
\begin{equation}\label{1}
\begin{split}
{h}_{t}^i=\nabla f_t^i(y_t^i)+\delta_tp_t^i, 
\end{split}
\end{equation}
where the vector $p_t^i$ is defined as
\begin{align}
\nonumber&p_{t}^i[k]=(\textbf{e}^k)^{\top}p_t^i
\\\nonumber=&\frac{(\textbf{e}^k)^{\top}\Big[\nabla^2f_t^i\left(y_t^i+[\beta_t^{i,k}]^{+}[\theta_t^{i,k}]^{+} \delta_t\textbf{e}^k\right)[\theta_t^{i,k}]^{+}\Big]\textbf{e}^k}{2}\\
+\!&\frac{(\textbf{e}^k)^{\top}\Big[\nabla^2f_t^i\left(y_t^i\!+\![\beta_t^{i,k}]^{\!-\!}[\theta_t^{i,k}]^{\!-\!} \delta_t\textbf{e}^k\right)[\theta_t^{i,k}]^{\!-\!}\Big]\textbf{e}^k}{2}\label{2}
\end{align}
with $[\theta_t^{i,k}]^{+},[\theta_t^{i,k}]^{-},[\beta_t^{i,k}]^{+},[\beta_t^{i,k}]^{-}\in [-1,1]$ for $k\in\{1,\ldots,d\}$. 
By Assumption 2,~4 and 5, from (\ref{1})--(\ref{2}), for $t\ge 1$, we can get 
\begin{equation}\label{3aa}
\begin{split}
\left\|{p}^i_t\right\|
\le \sqrt{d} L_H. 
\end{split}
\end{equation}
Define $L_h\triangleq L+\delta_0\sqrt{d} L_H$, then for $t\ge 0$, 
\begin{equation}\label{3a}
\begin{split}
\left\|{h}^i_t\right\|\le\left\|\nabla f_t^i(y_t^i)\right\|+\delta_t\left\|p_t^i\right\|\leq  L_h. 
\end{split}
\end{equation}
For $t\geqslant 0$, define $q^i_{t+1}=x^i_{t+1}-y^i_t$. By (\ref{sa4}), (\ref{Proj1}), (\ref{theta}), and (\ref{3a}), it follows that
\begin{align}
\nonumber\left\|q^i_{t+1}\right\|&=  \left\|\mathcal{P}_{X_{\delta_t}}\left(y^i_t-\alpha_t \frac{h^i_t}{ z_{t+1}[i] }\right)-\mathcal{P}_{X_{\delta_t}}\left(y^i_t\right)\right\| \\
& \le \frac{\alpha_t}{ z_{t+1}[i] }\left\|h_t^i\right\| \le \xi_1 {L}_{h} \alpha_t\label{4}
\end{align}
(\ref{2})--(\ref{4}) and similar lines of Theorem \ref{sThm1} yields Theorem \ref{Thm1}. {Note that by the definition of $\mathcal{P}_X(\cdot)$ and \eqref{3a},
\begin{align}
\nonumber\left\|q_{t+1}^i+\alpha_t\frac{h_t^i}{z_{t+1}[i]}\right\|&=\left\|x_{t+1}^i-\left(y_t^i-\alpha_t\frac{h_t^i}{z_{t+1}[i]}\right)\right\|
\\&\leq \alpha_t\frac{\left\|h_t^i\right\|}{z_{t+1}[i]}\leq \frac{L_h\alpha_t}{z_{t+1}[i]}, \label{10a}
\end{align}
{which combines (\ref{Proj2}) in Lemma \ref{Lem2.3} and \eqref{eq:pi=s} yields 
\begin{align}
\nonumber&\pi_{t+1}[i]\left\langle q_{t+1}^i+\alpha_t\frac{h_t^i}{z_{t+1}[i]}, x_{t+1}^i-{x}_t^*\right\rangle
\\=&\nonumber\pi_{t+1}[i]\left\langle q_{t+1}^i+\alpha_t\frac{h_t^i}{z_{t+1}[i]}, x_{t+1}^i-\tilde{x}_t^*+\tilde{x}_t^*-{x}_t^*\right\rangle
\\\leq& \frac{\pi_{t+1}[i]}{z_{t+1}[i]}L_h\alpha_t \left\|\tilde{x}_t^*-{x}_t^*\right\|\leq \frac{ML_h}{r} \alpha_t\delta_t
\label{11}
\end{align}
Similar to the proof in (\ref{s9})-(\ref{s13}) and replace $\tilde{x}_t^*$ by $x_t^*$, 
\begin{align}
&\nonumber\frac{1}{2}\left\|\bar{x}_{t\!+\!1}\!-\!{x}_{t\!+\!1}^*\right\|^2\!-\!\frac{1}{2}\left\|\bar{x}_t\!-\!{x}_t^*\right\|^2
\!\le \!2M\theta_t\!+\!2\xi_1L_h\alpha_t\beta_{t\!+1\!}
\\&+\!\frac{1}{2}\xi_1^2\alpha_t^2L_h^2\!+\!\alpha_t\sum_{i=1}^n\left\langle h_t^i, {x}_t^*\!-\!\bar{x}_t\right\rangle\!+\!\frac{nML_h}{r} \alpha_t\delta_t.\label{17}
\end{align}}
By (\ref{1}), we have 
\begin{align}
\nonumber& \sum_{i=1}^n \left\langle h_t^i, x_t^*-\!\bar{x}_t\right\rangle 
= \sum_{i=1}^n\left(\left\langle\nabla f_t^i\left(y_t^i\right), x_t^*-y_t^i\right\rangle\right.\\&\left.+\left\langle\nabla f_t^i\left(y_t^i\right), y_t^i-\bar{x}_t\right\rangle +\left\langle \delta_t p_t^i, x_t^*-\bar{x}_t\right\rangle\right)
 %\!=\!E(1)\!+\!E(2)\!+\!E(3), 
 \label{14}
\end{align}
%with 
%\begin{align}
%E_t(1)=& \sum_{i=1}^n \frac{\pi_{t+1}[i]}{z_{t+1}[i]}\left\langle\nabla f_t^i\left(y_t^i\right), x_t^*-y_t^i\right\rangle, \\
%E_t(2)=& \sum_{i=1}^n \frac{\pi_{t+1}[i]}{z_{t+1}[i]}\left\langle\nabla f_t^i\left(y_t^i\right), y_t^i-\bar{x}_t\right\rangle, 
%\end{align}
%and 
%\begin{align}
%E_t(3)=\sum_{i=1}^n \frac{\pi_{t+1}[i]}{z_{t+1}[i]}\left\langle c_t p_t, x_t^*-\bar{x}_t\right\rangle. 
%\end{align}
By a similar statement of (\ref{s13b}) and Schwarz inequality, Assumption 2, we can get  %for $E_t(2)$ 
\begin{align}
%&E_t(2)=
\sum_{i=1}^n \left\langle\nabla f_t^i\left(y_t^i\right), y_t^i-\bar{x}_t\right\rangle
%\\\le & \sum_{i=1}^n \frac{\pi_{t+1}[i]}{z_{t+1}[i]}\left\|\nabla f_t^i\left(y_t^i\right)\right\|\cdot \left\|y_t^i-\bar{x}_t\right\|
\le nL\beta_t\label{13bb}
\end{align}
By Schwarz inequality, (\ref{3aa}), Assumption 1, (\ref{theta}) and $\pi_t$ is stochastic, %for $E_t(3)$
we can obtain
\begin{align}
\sum_{i=1}^n \left\langle \delta_t p_t^i, x_t^*-\bar{x}_t\right\rangle\le 2Mn\sqrt{d}L_H\delta_t
\label{16}
\end{align}
By (\ref{14})--(\ref{16}) and the similar statements of (\ref{s13b})--(\ref{s17}), we have 
\begin{align}
\nonumber&\sum_{i=1}^n\langle h_t^i,x_t^*-\bar x_t\rangle
\le \sum_{i=1}^n\left[f_t^i(x_t^*)-f_t^i(\bar x_t)
-\frac\sigma 2\|\bar x_t-x_t^*\|^2\right]\\
&~~~+n(2L+2M\sigma)\beta_t
+2Mn\sqrt dL_H\delta_t.\label{17a}
\end{align}
By (\ref{3a})--(\ref{17a}) and the similar statements of %(\ref{s13b})--(\ref{s17}), 
(\ref{s23})--(\ref{s29}), we can get (\ref{18}) and  (\ref{18a}). \hfill$\square$

\section{Proof of Corollary~\ref{cor:static}}\label{Proof:static}
{For time-invariant objectives, $x_t^*=x^*$, $\theta_t=0$, and $\tilde x_t^*=(1-\delta_t/r)x^*$ for all $t \ge 0$. %define $f=\sum_{i=1}^n f^i$. 

For Algorithm~\ref{alg:rand}, define $\mathcal{F}_t\triangleq \sigma\{u_s^i:~0\le s\le t-1, i\in V\}$ for  $t\ge 1$, where $u_s^i$ is the random variable defined in Algorithm 1. Thus by Algorithm 1 we can get $x_t^i, y_t^i$ are $ \mathcal{F}_t-$measurable for $t\ge 1,~i\in V$, \eqref{aa2'} and the fresh uniform draw give
\[
\mathbb E[g_t^i\mid\mathcal F_t]
=\nabla\widehat f^i_{\delta_t}(y_t^i), 
\]
where $\widehat f^i_{\delta_t}(x)\triangleq\mathbb{E}[f^i(x+\delta_t v)]$ with $v$ uniformly distributed over $\mathbb{B}$ for $x\in X_{\delta_t}$. Taking conditional expectations in \eqref{s13} and applying the bounds in the proof of Lemma~\ref{sLem4.2} gives
\begin{align}
\notag &\tfrac12\mathbb E[\|\bar x_{t+1}-\tilde x_{t+1}^*\|^2\mid\mathcal F_t]\\
\notag &\le \tfrac12(1-n\sigma\alpha_t)\|\bar x_t-\tilde x_t^*\|^2
 -\alpha_t[f(\bar x_t)-f(x^*)]\\
\notag &\quad+\tfrac12\xi_1^2 L_g^2\alpha_t^2
 +2\xi_1 L_g\alpha_t\beta_{t+1}
 +n(2M\sigma+L+L_g)\alpha_t\beta_t\\
&\quad+C_4\alpha_t\delta_t+\tfrac{2M^2}{r}(\delta_t-\delta_{t+1}).
\label{static:rand}
\end{align}
Take $\sigma=0$ in the convex case. By Cauchy--Schwarz inequality, we have 
\[
\sum_{t=3}^{\infty}
 \left(\sum_{s=2}^{t-1}\alpha_{s-1}\lambda^{t-s-1}\right)^2
 \le \frac{1}{(1-\lambda)^2}\sum_{s=1}^{\infty}\alpha_s^2<\infty.
\]
Thus the definition  of $\beta_t$ in \eqref{s86} gives $\sum_{t=1}^\infty\beta_t^2<\infty$, and hence
\begin{align*}
\sum_{t=1}^\infty\alpha_t(\beta_t+\beta_{t+1})
&\le 2\left(\sum_{t=1}^\infty\alpha_t^2
             \sum_{t=1}^\infty\beta_t^2\right)^{1/2}<\infty,\\
\sum_{t=1}^\infty(\delta_t-\delta_{t+1})&=\delta_1 .
\end{align*}
Since $\sum_{t=1}^\infty\alpha_t\delta_t<\infty$, applying the Robbins--Siegmund almost-supermartingale convergence theorem (see, Theorem 3.1 in ~\cite{2010Ram}) to \eqref{static:rand} yields, almost surely,
\begin{align}
\notag &\{\|\bar x_t-\tilde x_t^*\|^2\}_{t\ge1}\ \text{converges},\\
&\sum_{t=1}^\infty\alpha_t[f(\bar x_t)-f(x^*)]<\infty .
\label{static:gap}
\end{align}
Note that
\[
\big|\|\bar x_t-x^*\|-\|\bar x_t-\tilde x_t^*\|\big|
\le \|\tilde x_t^*-x^*\|\le M\delta_t/r, 
\]
which combines Assumption 4 yields  $\|\bar x_t-x^*\|$ has a limit almost surely. Since $\sum_{t=1}^\infty\alpha_t=\infty$, the nonnegative gap in \eqref{static:gap} satisfies
\[
\liminf_{t\to\infty}[f(\bar x_t)-f(x^*)]=0
\quad\text{a.s.}
\]
By compactness of $X$, continuity of the fixed objective, and uniqueness of its minimizer, a subsequence of $\bar x_t$ therefore converges to $x^*$ almost surely. Thus the limit of $\|\bar x_t-x^*\|$ is zero. Theorem~\ref{sThm1} gives, for all $i\in V$,
\[
\|x_t^i-x^*\|\le\beta_t+\|\bar x_t-x^*\|\to0
\quad\text{a.s.}
\]

For Algorithm~\ref{alg:det}, by \eqref{17} and \eqref{17a} we can get %for every admissible sequence of graphs and attacks,
\begin{align}
\notag &\tfrac12\|\bar x_{t+1}-x^*\|^2
 \le\tfrac12(1-n\sigma\alpha_t)\|\bar x_t-x^*\|^2\\
\notag &\quad-\alpha_t[f(\bar x_t)-f(x^*)]
 +\tfrac12\xi_1^2L_h^2\alpha_t^2+2\xi_1L_h\alpha_t\beta_{t+1}\\
&\quad+n(2L+2M\sigma)\alpha_t\beta_t+C'_4\alpha_t\delta_t .
\label{static:det}
\end{align}
Taking $\sigma=0$ and summing \eqref{static:det} from $s=t$ to $T$ yields
\begin{align*}
&\tfrac12\|\bar x_{T+1}-x^*\|^2
 +\sum_{s=t}^T\alpha_s[f(\bar x_s)-f(x^*)]\\
&\le\tfrac12\|\bar x_t-x^*\|^2
 +O\!\left(\sum_{s=t}^\infty\alpha_s
       (\alpha_s+\beta_s+\beta_{s+1}+\delta_s)\right).
\end{align*}
By \eqref{86} and the same summability estimates, the tail tends to zero as $t\to\infty$. Taking $t=1$ and letting $T\to\infty$ yields
$\sum_{s=1}^{\infty}\alpha_s[f(\bar{x}_s)-f(x^*)]<\infty$.  Dropping the nonnegative sum on the left, taking the limit superior in $T$, and then the limit inferior in $t$ gives
\[
\limsup_{T\to\infty}\|\bar x_{T+1}-x^*\|^2
\le \liminf_{t\to\infty}\|\bar x_t-x^*\|^2 .
\]
The same compactness and uniqueness argument and Theorem~\ref{Thm1} now give, deterministically,
\[
\|x_t^i-x^*\|\le\beta_t+\|\bar x_t-x^*\|\to 0,\qquad i\in V, 
\]
which completes the proof of (i). 

If $f^i$ is $\sigma-$strongly convex for $i\in V$,  by $n\sigma$-strong convexity of $f$ and optimality of $x^*$ on $X$,
\[
 f(\bar x_t)-f(x^*)\ge\tfrac{n\sigma}{2}\|\bar x_t-x^*\|^2 .
\]
For Algorithm~\ref{alg:rand}, since $\|\tilde x_t^*-x^*\|\le M\delta_t/r$, we have
\begin{align*}
&\left|\tfrac12\|\bar x_t-x^*\|^2-\tfrac12\|\bar x_t-\tilde x_t^*\|^2\right|\\
&=\tfrac12\left|\langle 2\bar x_t-x^*-\tilde x_t^*,\,\tilde x_t^*-x^*\rangle\right|
 \le 2M^2\delta_t/r .
\end{align*}
Consequently,
\[
 f(\bar x_t)-f(x^*)\ge
 \tfrac{n\sigma}{2}\|\bar x_t-\tilde x_t^*\|^2-\tfrac{2n\sigma M^2}{r}\delta_t .
\]
Substituting this bound into \eqref{s8} gives
\begin{align*}
&\tfrac12\mathbb E\|\bar x_{t+1}-\tilde x_{t+1}^*\|^2
 \le(\tfrac12-n\sigma\alpha_t)\mathbb E\|\bar x_t-\tilde x_t^*\|^2\\
&\quad+\tfrac12\xi_1^2 L_g^2\alpha_t^2
 +2\xi_1L_g\alpha_t\beta_{t+1}
 +n(2M\sigma+L+L_g)\alpha_t\beta_t\\
&\quad+(C_4+2n\sigma M^2/r)\alpha_t\delta_t
 +\tfrac{2M^2}{r}(\delta_t-\delta_{t+1}).
\end{align*}
For Algorithm~\ref{alg:det}, substituting the strong-convexity bound directly into \eqref{static:det} gives
\begin{align*}
&\tfrac12\|\bar x_{t+1}-x^*\|^2
 \le(\tfrac12-n\sigma\alpha_t)\|\bar x_t-x^*\|^2\\
&\quad+\tfrac12\xi_1^2L_h^2\alpha_t^2+2\xi_1L_h\alpha_t\beta_{t+1}
 +n(2L+2M\sigma)\alpha_t\beta_t\\
&\quad+C'_4\alpha_t\delta_t .
\end{align*}
For the stated $\alpha_t$ and $\delta_t$, splitting the sum in \eqref{s86} and \eqref{86} at $s=\lfloor t/2\rfloor$ yields
\[
\sum_{s=2}^{t-1}\alpha_{s-1}\lambda^{t-s-1}
 =O(t\lambda^{t/2})+O(1/t)=O(1/t).
\]
Thus $\beta_t=O(1/t)$ for both algorithms. Using
\[
1-2n\sigma\alpha_t=\tfrac{t-1}{t+1},\qquad
\delta_t-\delta_{t+1}=\tfrac{C_\delta}{(t+1)(t+2)},
\]
the two recursions reduce to
\begin{align*}
\mathbb E\|\bar x_{t+1}-\tilde x_{t+1}^*\|^2
&\le\tfrac{t-1}{t+1}\mathbb E\|\bar x_t-\tilde x_t^*\|^2+O(t^{-2}),\\
\|\bar x_{t+1}-x^*\|^2
&\le\tfrac{t-1}{t+1}\|\bar x_t-x^*\|^2+O(t^{-2}).
\end{align*}
Multiplying by $t(t+1)$ and summing from $t=1$ to $T$ gives
\begin{align*}
T(T+1)\mathbb E\|\bar x_{T+1}-\tilde x_{T+1}^*\|^2&=O(T),\\
T(T+1)\|\bar x_{T+1}-x^*\|^2&=O(T).
\end{align*}
Dividing by $T(T+1)$ and using the consensus bounds, for each $i\in V$ we obtain, for Algorithm~\ref{alg:rand},
\begin{align*}
\mathbb E\|x_t^i-x^*\|^2
&\le 3\beta_t^2+3\mathbb E\|\bar x_t-\tilde x_t^*\|^2
       +3M^2\delta_t^2/r^2\\
&=O(t^{-2})+O(t^{-1})+O(t^{-2})=O(1/t),
\end{align*}
and, for Algorithm~\ref{alg:det},
\[
\|x_t^i-x^*\|^2\le 2\beta_t^2+2\|\bar x_t-x^*\|^2=O(1/t), 
\]
which completes the proof. \hfill$\square$}

\bibliographystyle{plain}        % Include this if you use bibtex 
\bibliography{ref}           % and a bib file to produce the 
                                 % bibliography (preferred). The
                                 % correct style is generated by
                                 % Elsevier at the time of printing.

%\begin{thebibliography}{99}     % Otherwise use the  
                                 % thebibliography environment.
                                 % Insert the full references here.
                                 % See a recent issue of Automatica 
                                 % for the style.
%  \bibitem[Heritage, 1992]{Heritage:92}
%     (1992) {\it The American Heritage. 
%     Dictionary of the American Language.}
%     Houghton Mifflin Company.
%  \bibitem[Able, 1956]{Abl:56}
%     B.~C.~Able (1956). Nucleic acid content of macroscope. 
%     {\it Nature 2}, 7--9. 
%  \bibitem[Able {\em et al.}, 1954]{AbTaRu:54}   
%     B.~C. Able, R.~A. Tagg, and M.~Rush (1954).
%     Enzyme-catalyzed cellular transanimations.
%     In A.~F.~Round, editor, 
%     {\it Advances in Enzymology Vol. 2} (125--247). 
%     New York, Academic Press.
%  \bibitem[R.~Keohane, 1958]{Keo:58}
%     R.~Keohane (1958).
%     {\it Power and Interdependence: 
%     World Politics in Transition.}
%     Boston, Little, Brown \& Co.
%  \bibitem[Powers, 1985]{Pow:85}
%     T.~Powers (1985).
%     Is there a way out?
%     {\it Harpers, June 1985}, 35--47.

%\end{thebibliography}

\end{document}